\input amssym.tex
\magnification 1200
\hsize = 14.5cm
\hoffset -0.5cm

\font\Bbb=msbm10
\def\BBB#1{\hbox{\Bbb#1}}
\font\Frak=eufm10
\def\frak#1{{\hbox{\Frak#1}}}
\font\smallFrak=eufm10 scaled 800
\def\smallfrak#1{{\hbox{\smallFrak#1}}}

\font\bol=cmbx10 scaled 1200

\hyphenation{Bor-cherds}
\hyphenation{pre-print}


\def\Ldg{1.1}
\def\Ldk{1.2}
\def\Ldd{1.3}

\def\kqu{2.1}
\def\dqu{2.2}
\def\compa{2.3}
\def\compb{2.4}
\def\TT{2.5}
\def\dvw{2.6}
\def\dovw{2.7}
\def\gvw{2.8}

\def\vir{3.1}
\def\comm{3.2}
\def\degr{3.3}
\def\omdeg{3.4}
\def\Bora{3.5}
\def\Borb{3.6}
\def\skeww{3.7}
\def\Ytens{3.8}
\def\omtens{3.9}
\def\vla{3.10}
\def\Y{3.11}
\def\LL{3.12}
\def\Lf{3.13}
\def\fg{3.14}
\def\LLz{3.15}
\def\Lfz{3.16}
\def\fgz{3.17}
\def\cvp{3.18}
\def\vhei{3.19}
\def\newVir{3.20}
\def\hvp{3.21}
\def\xyK{3.22}
\def\twi{3.23}
\def\xoe{3.24}
\def\Ye{3.25}
\def\Yw{3.26}

\def\korz{4.1}
\def\grz{4.2}
\def\dorz{4.3}
\def\tdp{4.4}
\def\tdo{4.5}
\def\tdadb{4.6}
\def\tdodb{4.7}
\def\tdodo{4.8}
\def\Rkakb{4.9}
\def\Rgka{4.10}
\def\Rgg{4.11}
\def\Rdjg{4.12}
\def\Rdog{4.13}
\def\Rdidj{4.14}
\def\Rdodj{4.15}
\def\Rdodo{4.16}
\def\diag{4.17}
\def\RS{4.18}
\def\koq{4.19}
\def\kpq{4.20}
\def\doq{4.21}
\def\tdpq{4.22}
\def\dpq{4.23}
\def\Eab{4.24}
\def\Yko{4.25}
\def\Ygm{4.26}
\def\Ykam{4.27}
\def\Ydam{4.28}
\def\Ydom{4.29}
\def\gamo{4.30}
\def\kbg{4.31}
\def\gpg{4.32}
\def\kdbg{4.33}
\def\Ykg{4.34}
\def\Ykkg{4.35}
\def\kbao{4.36}
\def\kdoao{4.37}
\def\kdbao{4.38}
\def\kaq{4.39}
\def\Ykkao{4.40}
\def\Ykko{4.41}
\def\koqq{4.42}
\def\koda{4.43}
\def\kbda{4.44}
\def\gda{4.45}
\def\doda{4.46}
\def\dbda{4.47}
\def\Epam{4.48}
\def\gE{4.49}
\def\doE{4.50}
\def\dbE{4.51}
\def\daq{4.52}
\def\damz{4.53}
\def\Eabo{4.54}
\def\Emab{4.55}
\def\qE{4.56}
\def\daqq{4.57}
\def\EE{4.58}

\def\hhh{5.1}


\def\catl{2.1}
\def\JM{{2.2}}
\def\RB{{2.3}}
\def\BB{{2.4}}
\def\BC{{2.5}}

\def\voa{3.1}
\def\inv{{3.2}}
\def\Cinv{{3.3}}
\def\RVL{{3.4}}
\def\tav{{3.5}}
\def\reln{{3.6}}
\def\RVA{{3.7}}
\def\FA{{3.8}}
\def\FAA{{3.9}}
\def\Rp{{3.10}}

\def\tvla{4.1}
\def\tops{{4.2}}
\def\TA{{4.3}}
\def\LAA{{4.4}}
\def\LAB{{4.5}}
\def\LAC{{4.6}}
\def\LAD{{4.7}}
\def\LAE{{4.8}}
\def\LAF{{4.9}}
\def\LAG{4.10}

\def\PAB{5.1}
\def\PA{{5.2}}
\def\TC{{5.3}}
\def\FB{{5.4}}
\def\sll{{5.5}}


\def\L{{{\cal L}}}
\def\D{{\cal D}}
\def\CC{{\cal C}}
\def\ZZ{{\cal Z}}
\def\U{{\cal U}}
\def\F{{\cal F}}
\def\B{{\cal B}}

\def\d{\partial}
\def\g{{\frak g}}
\def\f{{\frak f}}
\def\sg{{\smallfrak g}}
\def\sf{{\smallfrak f}}
\def\gl{{\it gl}}
\def\glN{{\it gl}_N}
\def\wgl{{\widehat {gl}_N}}
\def\slN{{{\it sl}_N}}
\def\sln{{{sl}_N}}
\def\wsl{{\widehat {sl}_N}}
\def\dg{{\dot \g}}
\def\sdg{{\dot \sg}}

\def\wdg{{\widehat \dg}}
\def\swdg{{\widehat \sdg}}
\def\td{\tilde d}

\def\df{{\dot \f}}

\def\tf{{\tilde \f}}

\def\om{\omega}
\def\ot{\otimes}

\def\div{{\rm div}}

\def\Hyp{{\mit Hyp}}
\def\hyp{{\mit Hyp}}

\def\VH{{V_{\Hyp}^+}}

\def\Hei{{{\cal H}{\mit ei}}}
\def\HV{{\cal VH}}
\def\R{{\cal R}}
\def\K{{\cal K}}

\def\C{\BBB C}
\def\o{{\bf 1}}
\def\Z{\BBB Z}

\def\r{{r}}
\def\m{{m}}
\def\s{{s}}
\def\u{{u}}
\def\vv{{v}}
\def\t{{t}}

\def\z{\left[ z_1^{-1} \delta \left( {z_2 \over z_1} \right) \right]}

\def\zd{\left[ z_1^{-1}  {\d \over \d z_2} 
\delta \left( {z_2 \over z_1} \right) \right]}

\def\zdd{\left[ z_1^{-1} \left( {\d \over \d z_2} \right)^2
\delta \left( {z_2 \over z_1} \right) \right]}

\def\zddd{\left[ z_1^{-1} \left( {\d \over \d z_2} \right)^3
\delta \left( {z_2 \over z_1} \right) \right]}

\def\dzb{{\d\over\d z_2}}
\def\dz{{\d\over\d z}}

\def\End{{\rm End}}
\def\Span{{\rm Span}}
\def\Id{{\rm Id}}
\def\Ind{{\rm Ind}}
\def\rank{{{\rm rank}\hbox{\hskip 0.08cm}}}

\def\Der{{{\rm Der}\hbox{\hskip 0.08cm}}}
\def\Ker{{\rm Ker \; }}
\def\Vir{{{\cal V}{\mit ir}}}

\def\deg{{{\rm deg}\hbox{\hskip 0.08cm}}}        
\def\char{{{\rm char}\hbox{\hskip 0.08cm}}}

\def\Inv{{inv}}

\

\centerline
{\bol A category of modules for the full toroidal Lie algebra.}

\

\centerline{
{\bf Yuly Billig}
\footnote{}{2000 Mathematics Subject Classification.
Primary 17B65, 17B69; Secondary 17B66.}
\footnote{}{Research supported by the  Natural Sciences and
 Engineering Research Council of Canada.}
}

\

\centerline
{\it To Robert Moody}


\

\

{\bf  Introduction.}

\

Toroidal Lie algebras are very natural multi-variable generalizations of 
affine Kac-Moody algebras. The theory of affine Lie algebras is rich and 
beautiful, having connections with diverse areas of mathematics and physics.
Toroidal Lie algebras are also proving themselves to be useful for the 
applications. Frenkel, Jing and Wang [FJW] used representations of toroidal Lie 
algebras to construct a new form of the  McKay correspondence.
Inami et al., studied toroidal symmetry in the context of a 4-dimensional 
conformal field theory [IKUX], [IKU]. There are also applications of toroidal
Lie algebras to soliton theory.  Using representations of the toroidal algebras
one can construct hierarchies of non-linear PDEs [B2], [ISW]. In particular, the toroidal
extension of the Korteweg-de Vries hierarchy contains the Bogoyavlensky's
equation, which is not in the classical KdV hierarchy [IT]. One can use the
vertex operator realizations to construct $n$-soliton solutions for the PDEs
in these hierarchies. We hope that further
development of the representation
theory of toroidal Lie algebras will help to find new applications of this
interesting class of algebras.

The construction of a toroidal Lie algebra is totally parallel to the well-known
construction of an (untwisted) affine Kac-Moody algebra [K1]. One starts with 
a finite-dimensional simple Lie algebra $\dg$ and considers Fourier polynomial maps from 
an $N+1$-dimensional torus into $\dg$. 
Setting $t_k = e^{ix_k}$, we may identify the algebra of Fourier polynomials
on a torus with the Laurent polynomial algebra 
$\R = \C[t_0^\pm, t_1^\pm, \ldots, t_N^\pm]$, and the
Lie algebra of the $\dg$-valued maps from a torus with the multi-loop algebra
$\C[t_0^\pm, t_1^\pm, \ldots, t_N^\pm] \ot \dg$. When $N = 0$, this yields
the usual loop algebra. 

Just as for the affine algebras, the next step is to build the universal central extension 
$(\R \ot \dg) \oplus \K$ of $\R \ot \dg$. 
However unlike the affine case, the center $\K$
is infinite-dimensional when $N \geq 1$. The infinite-dimensional center
makes this Lie algebra highly degenerate. One can show, for example, that in 
an irreducible bounded weight module, most of the center should act trivially.
To eliminate this degeneracy, we add the Lie algebra of vector fields on a torus,
$\D = \Der (\R)$ to $(\R \ot \dg) \oplus \K$. The resulting algebra,
$$\g = (\R \ot \dg) \oplus \K \oplus \D$$
is called the full toroidal Lie algebra (see Section 1 for details). 
The action of $\D$ on $\K$ is non-trivial,
making the center of the toroidal Lie algebra $\g$ finite-dimensional. This enlarged
algebra will have a much better representation theory.

 The most important class of modules for the affine Lie algebras are the highest weight
modules, and one would certainly want to construct their toroidal analogs.
The first problem that arises here is that one needs a triangular decomposition 
for the Lie algebra in order to introduce the notion of the highest weight module.
Toroidal Lie algebras are graded by $\Z^{N+1}$, and for $N > 0$, there is no canonical 
way of dividing this lattice into positive and negative parts. This difficulty is not 
present for the affine Lie algebras, which are graded by $\Z$, and for $\Z$ such a splitting 
is natural.

 One way to split $\Z^{N+1}$ is to cut it with a hyperplane that intersects with the lattice
only at zero. The corresponding class of the highest weight modules was studied by
Berman and Cox [BC], where it was found that the Verma modules constructed in this way 
will have infinite-dimensional weight spaces and do not produce any representations with
interesting realizations. Modules corresponding to other decompositions of the lattice
were studied in [DFP].

 An extremal way of dividing the lattice is to cut it with a hyperplane that
intersects $\Z^{N+1}$ at a sublattice of rank $N$. This approach was taken by Moody, Rao
and Yokonuma, who constructed a homogeneous vertex operator realization of the basic module for
the universal central extension of the multi-loop algebra [MRY]. In [EM], Rao and Moody
showed how to get a representation of a bigger algebra on the same space, adding a subalgebra
$\D^* = \mathop\oplus\limits_{p=1}^N \R {\d \over \d t_p}$ of the Lie algebra of vector fields.
One unanticipated development in [EM] was the appearance of a $\K$-valued 2-cocycle $\tau_1$ 
on the Lie algebra of vector fields, which is an abelian generalization of the Virasoro
cocycle.
 A principal realization for the basic module was given in [B1].
  
 Developing further these ideas, Larsson succeeded in constructing a wider class of 
representations for the toroidal Lie algebras [L]. He showed that the basic module
for the affine Lie algebra $\wdg$ may be replaced with an arbitrary highest weight module.
Larsson also discovered that affine $\wgl$-modules can be used as an ingredient in these constructions.
In Larsson's paper a combination of 2-cocycles $\tau_1$ and $\tau_2$ had appeared. 

 Berman and Billig developed a categorical approach to the representation theory
of toroidal Lie algebras, introducing the generalized Verma modules for toroidal
Lie algebras [BB]. They developed a theory of Lie algebras with polynomial multiplication,
which allowed them to prove that simple quotients of the generalized Verma modules always
have finite-dimensional weight spaces. Realizations of these irreducible quotients were obtained 
using a modified version of Larsson's construction. 

 A vertex algebra interpretation of these results was given by Berman, Billig and Szmigielski
in [BBS].

 Although the generalized Verma modules could be defined for the full toroidal algebras and
the result of [BB] concerning finite-dimensional weight spaces also holds in full generality,
there was one substantial difficulty unresolved in all of these publications. It was not known
how to construct realizations for the modules over the full toroidal algebra -- the piece
$\R {\d \over \d t_0}$ was always missing. This piece corresponds to the energy-momentum 
tensor in quantum field theories, so it is important to have it represented.

 A class of modules for the full toroidal Lie algebras was constructed in [B4] (unpublished),
and in the present paper we completely solve this problem.
 
 Extended affine Lie algebras (EALAs) are another important family of algebras closely related 
to the toroidal Lie algebras. The main feature of the extended affine Lie algebras is the existence
of a non-degenerate symmetric invariant bilinear form. Such a form does not exist on the full
toroidal algebra, but it can be defined on its subalgebra    
$(\R \ot \dg) \oplus \K \oplus \D_\div$, where $\D_\div$ is the Lie algebra
of divergence zero vector fields on a torus. 
 The results of the present paper make it possible to develop the representation theory
for the toroidal EALA using restriction from the full toroidal algebras [B4], [B6].

 It is shown in [ABFP] that most of the  extended affine Lie algebras may be realized as 
twisted toroidal EALAs. The representation theory of the twisted toroidal EALAs is
studied by Billig and Lau in [BL].

 In the present paper we define a rather natural category $\B_\chi$ of bounded $\g$-modules
with finite-dimensional weight spaces with the central character $\chi$. Our goal is to
study irreducible modules in this category. We show that every irreducible module is
characterized by its top $T$ -- the highest eigenspace for the operator 
$d_0 = t_0 {\d \over \d t_0}$. The space $T$ is a submodule with respect to the subalgebra
$\g_0$ consisting of elements of $\g$ of degree zero with respect to $t_0$. Following [BB],
we define a generalized Verma module $M(T)$ and its irreducible quotient $L(T)$. We show
(Theorem \BC) that every irreducible module in category $\B_\chi$ is isomorphic to $L(T)$
for some irreducible $\g_0$-module $T$ with finite-dimensional weight spaces. Using the 
results of [JM], [E] and [B5], we see that such $\g_0$-modules are precisely those considered
in [BB] -- they are multi-loop modules with respect to $\R_N \otimes \dg$ and tensor modules
with respect to $\Der (\R_N)$.

 Once we get a description of the tops $T$, we wish to completely determine the structure 
of the $\g$-modules $L(T)$. This is done by constructing the vertex operator realizations
of these modules. The crucial observation here is that the full toroidal Lie algebras are
vertex Lie algebras. This allows us to construct the universal enveloping VOA $V_\sg$. We
show (Proposition \tops) that for a particular top $T_0$, the irreducible module $L(T_0)$
is a factor-VOA of $V_\sg$. Using the methods developed in [BB], we study the kernel of the
projection $V_\sg \rightarrow L(T_0)$. This kernel gives us valuable information about
$L(T_0)$. Once we determine that a vector $v \in V_\sg$ belongs to the kernel, we apply
the state-field correspondence $Y$ and conclude that $Y_{L(T_0)} (v,z) = 0$. In this way we
derive important relations that hold in $L(T_0)$. We use these relations to define a toroidal 
VOA $V(T_0)$ as a tensor product of a sub-VOA $\VH$ of a lattice VOA and a VOA $V_\sf$
corresponding to the twisted Virasoro-affine Lie algebra with $\df = \dg \oplus \glN$.     

 In the previous papers on this subject, the vertex operator realization for the toroidal
modules had to be essentially guessed. The significant difference in the present approach is
that we are able to derive all the properties of the vertex operator realizations from inside,
using only the relations in the toroidal Lie algebra $\g$ and its universal enveloping
vertex algebra $V_\sg$.

 The VOA $V(T_0)$ controls the representation theory of $\g$. We show (Theorem \TC) that
every irreducible module $L(T)$ in category $B_\chi$ is a simple VOA module for $V(T_0)$
and can be constructed as a tensor product of a simple module $M_\Hyp^+ (\alpha)$ for the VOA
$\VH$
with an irreducible highest weight module $L_\sf$ for the twisted Virasoro-affine algebra
$\f$. For a generic level $c$, we can further factor $L_\sf$ in a tensor product and get 
the following decomposition for $L(T)$:
$$ L(T) \cong M_\Hyp^+ (\alpha) \otimes L_\swdg \otimes L_\wsl \otimes L_\Hei \otimes L_\Vir ,$$
where the last four factors are certain irreducible highest weight modules
for the affine algebras $\wdg, \wsl$, the infinite-dimensional Heisenberg algebra and the
Virasoro algebra. In this way we reduce the representation theory of toroidal Lie algebras
to the representation theory of affine, Heisenberg and Virasoro algebras.
 Whenever explicit realizations are available for the components in the tensor product 
decomposition above, we get a realization for the irreducible module for the full toroidal Lie
algebra.

This leads to the following open problem: 
while the explicit expressions for the characters of irreducible modules may be known, 
there is no Weyl-type character formula for the toroidal Lie algebras. Obtaining such a formula 
may yield interesting number-theoretic identities. 

 The structure of the paper is the following. In Section 1 we review the construction of the 
toroidal Lie algebras. In Section 2 we introduce a category $\B_\chi$ of $\g$-modules and show
that every irreducible module in $\B_\chi$ is characterized by its top $T$. We also describe
the structure of the top $T$. In Section 3 we recall the definition of the vertex operator algebra
and the construction of the universal enveloping vertex algebra of a vertex Lie algebra.
We introduce twisted Virasoro-affine Lie algebras, show that these algebras are in fact vertex Lie 
algebras and describe the structure of the corresponding enveloping vertex algebras and their 
simple modules. At the end of Section 3 we describe the hyperbolic lattice VOA $V_\Hyp$ and its 
sub-VOA $\VH$. In Section 4 we show that the full toroidal Lie algebra is vertex Lie algebra
and define its enveloping VOA $V_\sg$. Next we establish a series of relations that hold in the
simple quotient $L(T_0)$ and use them to decompose $L(T_0)$ into a tensor product of two VOAs,
$\VH$ and $L_\sf (\gamma_0)$.
In the final Section 5 we prove that a slightly bigger VOA $V(T_0) = \VH \otimes V_\sf (\gamma_0)$ 
also admits the structure of a module over the full toroidal algebra $\g$.
We show that every irreducible $\g$-module in category $\B_\chi$ is
a simple VOA module for $V(T_0)$, which then allows us to obtain a complete description
of these irreducible $\g$-modules.

\

{\bf Acknowledgements:} I thank Stephen Berman for the stimulating discussions
and encouragement. I have greatly benefited from Chongying Dong's lectures on vertex operator 
algebras given at the Fields Institute. 

\

{\bf 1. Toroidal Lie algebras.}

\

Toroidal Lie algebras are the natural multi-variable generalizations of
affine Lie algebras. 
In this review of the toroidal Lie algebras we follow the work [BB].
Let $\dg$ be a simple finite-dimensional Lie algebra
over $\C$ with a non-degenerate invariant bilinear form $( \cdot | \cdot )$
and let $N \geq 1$ be an integer. 
We consider the Lie algebra $\R \ot \dg$
of maps from an $N+1$ dimensional torus into $\dg$, where  
$\R = \C [t_0^\pm, t_1^\pm, \ldots, t_N^\pm]$ is the algebra of Fourier polynomials
on a torus. The universal central extension
of this Lie algebra may be described by means of the following construction which is
due to Kassel [Kas]. Let $\Omega^1_\R$ be the space of $1$-forms on a torus:
$\Omega^1_\R = \mathop\oplus\limits_{p = 0}^N \R dt_p$. We will choose the
forms $\{ k_p = t_p^{-1} dt_p | p = 0, \ldots, N \}$ as a basis of this
free $\R$ module. There is a natural map $d$ from the space of functions
$\R$ into $\Omega^1_\R$: $d(f) = \sum\limits_{p = 0}^N {\d f \over \d t_p} dt_p
= \sum\limits_{p = 0}^N t_p {\d f \over \d t_p} k_p$. The center $\K$ for the
universal central extension $(\R \ot \dg) \oplus \K$ of $\R \ot \dg$
is realized as
$$ \K = \Omega^1_\R / d(\R) , $$
and the Lie bracket is given by the formula
$$[f_1(t) g_1, f_2(t) g_2] = f_1(t) f_2(t) [g_1, g_2] + (g_1| g_2) f_2 d(f_1).$$
Here and in the rest of the paper we will denote elements of $\K$ by the
same symbols as elements of $\Omega^1_\R$, keeping in mind the canonical
projection $\Omega^1_\R \rightarrow \Omega^1_\R / d(\R)$.

Next we add to $(\R \ot \dg) \oplus \K$ the algebra $\D$ of 
vector fields on the torus
$$\D = \mathop\oplus\limits_{p=0}^N \R d_p,$$
where $ d_p = t_p {\d \over \d t_p}$.
We will use the multi-index notation writing
$\t^\r =  t_0^{r_0} t_1^{r_1} \ldots t_N^{r_N}$ for 
$\r = (r_0, r_1, \ldots, r_N)$, etc.

The natural action of $\D$ on $\R \ot \dg$ 
$$[\t^\r d_a, \t^\m g] = m_a \t^{\r+\m} g \eqno{(\Ldg)}$$
uniquely extends to the action on the universal central extension
$(\R \ot \dg) \oplus \K$ by
$$[\t^\r d_a, \t^\m k_b] = m_a \t^{\r+\m} k_b + \delta_{ab}
\sum\limits_{p=0}^N r_p \t^{\r+\m} k_p . \eqno{(\Ldk)}$$
This corresponds to the Lie derivative action of the vector fields on 1-forms.

It turns out that there is still an extra degree of freedom in defining
the Lie algebra structure on $(\R \ot \dg) \oplus \K \oplus \D$.
The Lie bracket on $\D$ may be twisted with a $\K$-valued 2-cocycle:
$$[\t^\r d_a, \t^\m d_b] = m_a \t^{\r+\m} d_b - r_b \t^{\r+\m} d_a
+ \tau(\t^\r d_a, \t^\m d_b) . \eqno{(\Ldd)}$$
In order to compute the second cohomology space $H^2 (\D,\K)$,
one could use the Gelfand-Fuks cohomology theory [F], [T].
Unfortunately this theory only allows one to do the computations in the
$C^\infty$ setup, i.e., when $\R$ is replaced with the algebra of infinitely differentiable 
functions on a torus, and not for the algebra of Fourier polynomials
that we consider here. For the $C^\infty$ situation the calculation of
$H^2_{C^\infty} (\D,\K)$ has been carried out in [BN]. For the $(N+1)$-dimensional torus
with $N+1 \geq 2$, the dimension of the second cohomology space is
$$\dim H^2_{C^\infty} (\D, \K) = 2 + \pmatrix{ N+1 \cr 3 \cr}, $$
and the following cocycles form the basis of this space:     
$$\tau_1 (\t^\r d_a, \t^\m d_b) = m_a r_b \sum\limits_{p=0}^N
m_p \t^{\r+\m} k_p ,$$
$$\tau_2 (\t^\r d_a, \t^\m d_b) = r_a m_b \sum\limits_{p=0}^N
m_p \t^{\r+\m} k_p ,$$
together with a family $\{ \eta_{abc} | 0 \leq a < b < c \leq N \}$,
where the cocycle $\eta_{abc}$ is defined by the conditions that
$$ \eta_{abc} (t^r d_{\sigma(a)}, t^m d_{\sigma(b)} ) = 
(-1)^\sigma t^{r+m} k_{\sigma(c)},$$
for any permutation $\sigma: \{ a, b, c \} \rightarrow \{ a, b, c \}$
and $\eta_{abc} ( t^r d_i, t^m d_j ) = 0$ if $i=j$ or $\{ i, j \} \not\subset
\{ a, b, c \}$.
It is clear that $H^2 (\D,\K)$ in the algebraic setup contains
the space $H^2_{C^\infty} (\D, \K)$. After twisting with a cocycle
$\eta_{abc}$, the vector fields $d_a$ and $d_b$ no longer commute. For this
reason we will be considering only the cocycles $\tau_1$ and $\tau_2$ in this paper.

We will write $\tau = \mu \tau_1 + \nu \tau_2$. 
The resulting algebra (or rather a family of algebras) is called the
full toroidal Lie algebra 
$$\g = \g(\mu,\nu) = (\R \ot \dg) \oplus \K \oplus \D.$$

Note that after adding the algebra of derivations $\D$, the center $\ZZ$
of the toroidal Lie $\g$ becomes finite-dimensional with the basis 
$\{ k_0, k_1, \ldots, k_N \}$. This can be seen from the action
(\Ldk) of $\D$ on $\K$,  which is non-trivial.

The study representation theory of toroidal Lie algebras has begun in [MRY]
and [EM], with further developments in [B1], [L], [BB], [BBS]. In all of these
papers there was one common difficulty that has not been resolved --
the representations constructed there were not for the full toroidal
algebra $\g$, but only for a subalgebra
$$ \g^* = (\R \otimes \dg) \oplus \K \oplus 
\left( \mathop\oplus\limits_{p=1}^N \R d_p \right),$$
where the piece $\R d_0$ that corresponds to the toroidal energy-momentum 
tensor was missing. 
This left the theory in a somewhat incomplete form, and the goal of the present
paper is to study representations for the full toroidal Lie algebra.  

\

\

{\bf 2. A category of bounded modules for toroidal Lie algebras.}

\

 In this section we will introduce a category
of bounded modules for toroidal Lie algebras that could be regarded as analogs
of the highest weight modules for affine Kac-Moody algebras. The difference
from the affine case is that the highest weight space is not 1-dimensional, but 
rather forms a multi-loop module for a smaller toroidal subalgebra.

These bounded modules are quite promising from the point of view of applications.
In [B2] a module of this type was used to construct a toroidal extension of
the Korteweg-de Vries hierarchy. 
The vertex operator realizations of the toroidal modules allow one to
construct soliton solutions to these non-linear PDEs.

The variable $t_0$ will play a special role in our construction. From the 
physics perspective, it may be interpreted as time, whereas $t_1, \ldots t_N$
are the space variables. 

 The algebra $\g$ has a $\Z^{N+1}$-grading by the eigenvalues of the adjoint
action of $d_0, d_1, \ldots d_N$. 
We will denote by $\{ \epsilon_0, \ldots, \epsilon_N \}$ the standard basis of
$\Z^{N+1}$.
We will be also considering its $\Z$-grading
just with respect to the action of $d_0$:
$$ \g = \mathop\oplus\limits_{n\in\Z} \g_n ,$$
and define subalgebras $\g_\pm = \mathop\oplus\limits_{n \gtrless 0} \g_n$,
which yields the decomposition $\g = \g_- \oplus \g_0 \oplus \g_+$.

We recall that the center $\ZZ$ of $\g$ is $N+1$-dimensional: 
$\ZZ = \Span \left< k_0, k_1, \ldots, k_N \right>$.
Clearly, in any irreducible weight module with finite-dimensional weight spaces, 
these central elements will act as multiplications by scalars. Let us fix a 
non-zero central character $\chi: \ZZ \rightarrow \C$ and define a category of bounded
modules with central character $\chi$ for the toroidal Lie algebra. 

{\bf Definition.}
A category  $\B_\chi$ of bounded modules for the toroidal Lie algebra is a 
category whose objects are $\g$-modules $B$ satisfying the following axioms:

(B1) $B$ has a weight decomposition with respect to the subalgebra 
$\left< d_0, d_1, \ldots, d_N \right>$ :
$$ B = \mathop\oplus\limits_{m\in\C^{N+1}} B_m ,$$
where 
$B_m = \left\{ v \in B \, | \, d_j (v) = m_j v, \; j =0,\ldots,N \right\};$

(B2) All weight spaces $B_m$ are finite-dimensional;

(B3) The action of the center $\ZZ$ on $B$ is given by the central character $\chi$:
$k_j v = \chi(k_j) v$ for all $v\in B, j=0,\ldots, N$;

(B4) Real parts of eigenvalues of $d_0$ on $B$ are bounded from above.

\

The physical meaning of the last axiom is that the spectrum of the energy operator
$E = -d_0$ has a lower bound, i.e., there exist states of the lowest energy.

The goal of this paper is to describe irreducible modules in category $\B_\chi$
and find their characters.

 First of all we are going to show that in order for $\B_\chi$ to be non-trivial,
the central character must satisfy the condition $\chi (k_1) = 0, \ldots,  
\chi(k_N) = 0.$ Clearly, $\chi$ must vanish on an $N$-dimensional subspace in $\ZZ$, 
and it turns out that this $N$-dimensional nullspace must be ``aligned'' with the
choice of the operator $d_0$ in axiom (B4).

{\bf Lemma \catl.} 
{\it
Suppose that $\B_\chi$ is a non-trivial category. 
Then  $\chi (k_j) = 0$ for all $j=1, \ldots, N$.    
}

{\it Proof.} Let $B$ be a non-zero module in $\B_\chi$. Let us reason by 
contradiction and assume that $\chi(k_j) = c_j \neq 0$ for some $j, \quad 1 \leq j \leq N$.

 Since the spectrum of $d_0$ is bounded, we can find a weight space $B_m$ 
such that $m_0+1$ is not an eigenvalue of $d_0$ on $B$. Let $v$ be a non-zero 
vector in $B_m$, and consider the following family of vectors:
$$ (t_0^{-1} t_j^{-n} k_j) (t_0^{-1} t_j^n k_j) v, \; \; \; n = 1, 2, \ldots$$
Clearly, all these vectors belong to the same weight space $B_{m-2\epsilon_0}$,
and we claim that they are all linearly independent. Indeed, suppose
$$\sum_{n > 0} a_n (t_0^{-1} t_j^{-n} k_j) (t_0^{-1} t_j^n k_j) v = 0 .$$
Since
$$ d_0 (t_0 t_j^r d_0) v = (m_0 +1) (t_0 t_j^r d_0) v, $$
and $m_0+1$ is not an eigenvalue of $d_0$ on $B$, we conclude that 
$(t_0 t_j^r d_0) v = 0$ for all $r\in\Z$. We also note that 
$$ [t_0 t_j^r d_0, t_0^{-1} t_j^s k_j ] = - t_j^{r+s} k_j = -\delta_{r,-s} k_j .$$
Taking these two facts into account, we get that for $r > 0$,
$$ 0 = (t_0 t_j^r d_0) (t_0 t_j^{-r} d_0) 
\sum_{n > 0} a_n (t_0^{-1} t_j^{-n} k_j) (t_0^{-1} t_j^n k_j) v = a_r c_j^2 v .$$
 Since $c_j \neq 0$, we conclude that $a_r = 0$ for all $r > 0$. 
Thus the vectors $\left\{ (t_0^{-1} t_j^{-n} k_j)
(t_0^{-1} t_j^n k_j) v \right\}$ with $n > 0$
are linearly independent, which contradicts (B2). 
This proves that $\chi (k_1) = 0, \ldots, \chi(k_N) = 0.$

\

 For the rest of this paper we fix a non-zero constant $c\in \C$ and let
$\chi = (c, 0, \ldots, 0)$.
From now on, the multivariable $t$ will not include $t_0$, that is
$t^r$ will stand for $t_1^{r_1} \ldots t_N^{r_N}$, etc.

 Consider an irreducible module $L$ in category $\B_\chi$. It is clear that the eigenvalues
of $d_0$ on $L$ belong to a single $\Z$-coset in $\C$. Let $d$ be the eigenvalue 
of $d_0$ with the highest real part, and let $T$ be the corresponding 
eigenspace.   
    
 Obviously, $T$ is a $\g_0$-module and $\g_+ T = 0$. 
It is easy to see that irreducibility of $L$ implies the irreducibility of 
$T$ as a $\g_0$-module. 
We will call the subspace
$T$ the {\it top} of $L$. Next we are going to describe the structure of $T$.
We will be using a result of [JM] for this.

{\bf Theorem {\JM} ([JM]).} 
{\it
Suppose $\chi(k_0) = c \neq 0, \; \chi(k_1) = 0, \ldots, \chi(k_N) = 0$.
Let $L$ be an irreducible module in category $\B_\chi$ with the top $T$.
Then 
$$ T \cong \C [q_1^\pm, \ldots q_N^\pm ] \otimes U,$$
where $U$ is a finite-dimensional space, and the action of $\g_0$ on $T$
satisfies
$$(t^\r k_0) (q^\m \otimes u) = c q^{\m+\r} \otimes u,
\quad\quad (t^\r k_j) (q^\m \otimes u) = 0, \eqno{(\kqu)}$$
$$ d_0 (q^m \otimes u) = d q^m \otimes u, \quad
d_j (q^\m \otimes u) = (m_j + \alpha_j) q^\m \otimes u, 
\quad u\in U, j=1,\ldots,N, \eqno{(\dqu)}$$
for some fixed $\alpha = (\alpha_1,\ldots,\alpha_N) \in \C^N, d\in\C$. 
}

If we take the quotient of $\g_0$ by the ideal  
$J = \Span\left\{ t^\r k_j | \r\in\Z^N, j = 1, \ldots, N \right\}$, 
which annihilates $T$,
then we will get a semi-direct product of the Lie algebra of vector fields
$\D_N = \Der \C[t_1^\pm,\ldots,t_N^\pm] $ on $N$-dimensional torus with 
a multi-loop algebra:
$$\g_0 / J \cong \D_N \ltimes \C[t_1^\pm,\ldots,t_N^\pm] \otimes 
\left( \dg \oplus \C d_0 \oplus \C k_0 \right) .$$

Since ${1\over c} (t^\r k_0)$ acts on $T$ as multiplication by $q^r$, we can derive
from (\Ldk) the following compatibility relations between the action of $\g_0$ and 
the operators of multiplication by $q^r$:
$$(t^\s d_j) q^\r - q^\r (t^\s d_j) = r_j q^{s+r} , \eqno{(\compa)}$$
$$ (t^\s d_0) q^r = q^\r (t^\s d_0), \quad\quad 
(t^\s g) q^r  = q^\r (t^\s g), \quad
 g\in\dg. \eqno{(\compb)}$$

 Eswara Rao [E] classified irreducible $\D_N$-modules with a compatible action
of the algebra of Laurent polynomials, proving that any such module is a tensor
module. We will use a version of this result for the semidirect 
product of $\D_N$ with a multi-loop algebra given in [B5], Theorem 4(c):

{\bf Theorem {\RB} ([E], [B5]).} 
{\it Let $\alpha \in \C^N, c, d \in \C$, $c\neq 0$. 
Let $T$ be an irreducible $\g_0$-module satisfying
the conclusion of Theorem \JM, as well as $(\compa),(\compb)$. 
Then there exist  a finite-dimensional irreducible $\dg$-module 
$V$ and a finite-dimensional irreducible $\gl_N$-module $W$, such that 
$$ T \cong \C [q_1^\pm, \ldots q_N^\pm ] \otimes V \otimes W , \eqno{(\TT)}$$    
and the action of $\g_0$ on  $T$ is given by (\kqu) and
$$(t^\r d_j) (q^\m \otimes v \otimes w) = 
(m_j + \alpha_j) q^{\m+\r} \otimes v \otimes w + 
\sum_{p=1}^N r_p q^{\m+\r} \otimes v \otimes E_{pj}w, 
\quad j = 1, \ldots, N, \eqno{(\dvw)}$$
$$ (t^r d_0) (q^\m \otimes v \otimes w) = d q^{m+r} \otimes v \otimes w,
\eqno{(\dovw)}$$
$$(t^r g) (q^\m \otimes v \otimes w) = q^{\m+\r} \otimes g v \otimes w , \quad g\in\dg . \eqno{(\gvw)}$$
}

Here in (\dvw) $E_{pj}$ denotes a matrix with 1 in position $(p,j)$, and zeros elsewhere.

Combining these two theorems, we conclude  
that an irreducible module in category $\B_\chi$ yields the following 
data -- a finite-dimensional irreducible $\dg$-module $V$, 
a finite-dimensional irreducible $\glN$-module $W$, 
a constant $d\in \C$ and $\alpha\in\C^N$. 
It is easy to see that the choice of $\alpha$ is not canonical and $\alpha$ can be changed 
to any value in the coset $\alpha+\Z^N$ by choosing a different weight space for the
generators of $T$ as a free $\C[q_1^\pm,\ldots,q_N^\pm]$-module.
An irreducible $\glN$-module $W$ is determined by the action of $\slN$ and a scalar
$h$, by which the identity matrix acts on $W$. 

 It was shown in [BB] that for any $\g_0$-module $T$, corresponding to
the data $(V,W,h,d,\alpha)$ as above, there exists
an irreducible module in $\B_\chi$ with $T$ as a top. Let us review this
construction.


First we let $\g_+$ act on $T$
trivially, and define the generalized Verma module as the induced module 
$$ M(T) = \Ind_{\sg_0\oplus\sg_+}^\sg (T) .$$
Note that the module $M(T)$ does not belong to category $\B_\chi$ since its weight spaces
lying below $T$ are infinite-dimensional. Nonetheless, the following result holds:

{\bf Theorem {\BB} ([BB]).} 
{\it
(a) The $\g$-module $M(T)$ has a unique maximal submodule $M^{rad}$.

(b) The factor-module $L(T) = M(T)/M^{rad}$ is an irreducible $\g$-module. 

(c) All weight spaces of $L(T)$ are finite-dimensional, and $L(T)$ belongs to the category
$\B_\chi$.
}

 \

 Summarizing, we get the following

{\bf Theorem \BC.} 
{\it (a) Let $\chi$ be a non-zero central character $\chi: \ZZ \rightarrow \C$.
A category $\B_\chi$ is non-trivial if and only if
$\chi(k_0) = c, \; \chi(k_1)= 0, \ldots, \chi(k_N)=0$ for some non-zero $c\in\C$.

 Let now $\chi = (c, 0, \ldots, 0)$ with $c\neq 0$.

(b) Irreducible $\g$-modules in category $\B_\chi$ are in 1-1 correspondence with 
the data 
\break
$(V, W, h, d, \alpha)$, where $V$ is a finite-dimensional irreducible 
$\dg$-module,
$W$ is a finite-dimensional irreducible $\slN$-module, $\alpha\in\C^N/\Z^N$, $h,d\in\C$.

(c) Every irreducible module in category $\B_\chi$ is isomorphic to $L(T)$ where 
$$T = \C [q_1^\pm, \ldots q_N^\pm ] \otimes V \otimes W$$
with the action of $\g_0$ given by (\kqu) and (\dvw)-(\gvw).  
}
  
{\it Proof.}  Part (a) has been already proved in Lemma {\catl}. Let us prove part (c).
As we have seen above, an irreducible module $L$ in category $\B_\chi$ has a top $T$,
the structure of which is described by Theorem \RB.
Thus $L$ is a factor-module of $M(T)$. However $M(T)$ has a unique
irreducible factor, which is isomorphic to $L(T)$. This proves $L \cong L(T)$.
Part (b) follows from (c) and Theorem \BB.

\

Our main goal will be to completely determine the structure of the irreducible modules $L(T)$,
and in particular, to find their characters. This will be done using the theory of vertex 
operator algebras (VOAs). We will show that for a particular choice of the data
$(V, W, h, d, \alpha)$, namely, with $V$ and $W$ being trivial 1-dimensional modules for
$\dg$ and $\slN$,
$\alpha = 0$, $h = N\nu c$ and $d= {1\over 2} (\mu + \nu) c$, yielding the top
$$ T_0 = \C [q_1^\pm, \ldots, q_N^\pm] ,$$
the module $L(T_0)$ is a vertex operator algebra, while all irreducible modules $L(T)$
are VOA-modules for a slightly bigger VOA $V(T_0)$. Once we determine the structure of 
$V(T_0)$ as a VOA, we will immediately get the structure of all the modules $L(T)$ using 
the principle of preservation of identities in the VOA theory.  

\

\

{\bf 3. Vertex operator algebras and vertex Lie algebras.}

\

{\bf 3.1. Definitions and properties of a VOA.}

Let us recall the basic notions of the theory of the vertex operator algebras.
Here we are following [K2] and [Li].

{\bf Definition.} 
{
A vertex algebra is a vector space
$V$ with a distinguished vector
$\o$ (vacuum vector) in $V$, 
an operator $D$ (infinitesimal translation) on the space $V$, and a linear map $Y$ (state-field
correspondence)
$$\eqalign{
Y(\cdot,z): \quad V &\rightarrow (\End V)[[z,z^{-1}]], \cr
a &\mapsto Y(a,z) = \sum\limits_{n\in\Z} a_{(n)} z^{-n-1} 
\quad (\hbox{\rm where \ } a_{(n)} \in \End V), \cr} $$
such that the following axioms hold:

\noindent
(V1) For any $a,b\in V, \quad a_{(n)} b = 0 $ for $n$ sufficiently large;

\noindent
(V2) $[D, Y(a,z)] = Y(D(a), z) = {d \over dz} Y(a,z)$ for any $a \in V$;

\noindent
(V3) $Y(\o,z) = \Id_V$;

\noindent
(V4) $Y(a,z) \o \in V [[z]]$ and $Y(a,z)\o |_{z=0} = a$ for any $a \in V$
\ (self-replication);

\noindent
(V5) For any $a, b \in V$, the fields $Y(a,z)$ and $Y(b,z)$ are mutually local, that is, 
$$ (z-w)^n \left[ Y(a,z), Y(b,w) \right] = 0, \quad \hbox{\rm for \ } 
 n  \hbox{\rm \ sufficiently large} .$$


A vertex algebra $V$ is called a vertex operator algebra (VOA) if, in addition, 
$V$ contains a vector $\omega$ (Virasoro element) such that

\noindent
(V6) The components $L(n) = \omega_{(n+1)}$ of the field
$$ Y(\omega,z) = \sum\limits_{n\in\Z} \omega_{(n)} z^{-n-1} 
= \sum\limits_{n\in\Z} L(n) z^{-n-2} $$
satisfy the Virasoro algebra relations:
$$  [ L(n) , L(m) ] = (n-m) L(n+m) + \delta_{n,-m} {n^3 - n \over 12} 
(\rank V) \Id, \quad \hbox{\rm where \ } \rank V \in \C;  
\eqno{(\vir)}$$

\noindent
(V7) $D = L(-1)$;

\noindent
(V8) $V$ is graded by the eigenvalues of $L(0)$:
$V = \mathop\oplus\limits_{n\in \Z} V_n$ with $L(0) \big|_{V_n} = n \Id$.
}

This completes the definition of a VOA.

As a consequence of the axioms of the vertex algebra  
we have the  following important commutator formula:
$$\left[ Y(a,z_1), Y(b,z_2) \right] =
\sum_{n \geq 0} {1 \over n!}  Y(a_{(n)} b, z_2)
\left[ z_1^{-1} \left( {\d \over \d z_2} \right)^n
\delta \left( {z_2 \over z_1} \right) \right] . \eqno{(\comm)}$$
 As usual, the delta function is
$$ \delta(z) = \sum_{n\in\Z} z^n .$$
By (V1), the sum in the right hand side of the commutator formula
is actually finite.

All the vertex operator algebras 
that appear in this paper have the gradings by 
non-negative integers (degree): $V = \mathop\oplus\limits_{n=0}^\infty V_n$.
In this case the sum in the right hand side of the commutator 
formula (\comm) runs from $n=0$ to $n = \deg(a) + \deg(b) -1$,
because 
$$\deg(a_{(n)} b) = \deg(a) + \deg(b) -n - 1 ,\eqno{(\degr)}$$ 
and the elements of negative degree vanish.

It follows from (V7) and (V8) that 
$$ \omega_{(0)} a = D(a), \quad
\omega_{(1)} a = \deg (a) a \quad \quad \hbox{\rm for } a
\hbox{\rm \  homogeneous}. \eqno{(\omdeg)} $$


Another consequence of the axioms of a vertex algebra is the Borcherds'
identity:
$$\sum\limits_{j\geq 0} \pmatrix{m \cr j \cr}
(a_{(k+j)} b)_{(m+n-j)} c $$
$$= \sum\limits_{j\geq 0} 
(-1)^{k+j+1} \pmatrix{k \cr j \cr}
b_{(n+k-j)} a_{(m+j)} c +
\sum\limits_{j\geq 0} 
(-1)^j \pmatrix{k \cr j \cr}
a_{(m+k-j)} b_{(n+j)} c, \quad \quad  k,m,n \in \Z. \eqno{(\Bora)} $$ 
We will particularly need its special case when $m=0$: 

 
$$ (a_{(k)} b)_{(n)} c = 
\sum\limits_{j\geq 0} (-1)^{k+j+1} \pmatrix{ k \cr j \cr}
b_{(n+k-j)} a_{(j)} c 
+ \sum\limits_{j\geq 0}  (-1)^j \pmatrix{k \cr j \cr}
a_{(k-j)} b_{(n+j)} c, \quad \quad  k,n \in \Z. \eqno{(\Borb)} $$


The last formula that we quote here is the skew-symmetry identity:
$$ a_{(n)} b = \sum_{j\geq 0} (-1)^{n+j+1} {1\over j!}
D^j (b_{(n+j)} a) . \eqno{(\skeww)}$$ 

\

{\bf 3.2. Tensor products of VOAs.}

Let us review here
the definition of the tensor product of two VOAs 
$\left( V^{\prime}, Y^{\prime}, \omega^{\prime}, \o \right)$
and $\left( V^{\prime\prime}, Y^{\prime\prime}, 
\omega^{\prime\prime}, \o \right)$
(the case
of an arbitrary number of factors is a trivial generalization).
The tensor product space $V = V^{\prime} \ot V^{\prime\prime}$ 
has the VOA structure under
$$ Y(a\ot b, z) = Y^{\prime} (a,z) \ot Y^{\prime\prime} (b,z), 
\eqno{(\Ytens)}$$
$$ \omega = \omega^{\prime} \ot \o 
+ \o \ot \omega^{\prime\prime},
\eqno{(\omtens)}$$
and $\o = \o \ot \o$ being the identity element.

It follows from (\omtens) that the rank of $V$ (see V6) is the sum of the ranks of the tensor factors. 

%
%
%
%
%

\

{\bf 3.3. Vertex Lie algebras.}

 An important source of the vertex algebras is provided by the vertex Lie algebras.
In presenting this construction we will be following [DLM] (see also [P], [R],
[K2], [FKRW]).

Let $\L$ be a Lie algebra with the basis 
$\{ u(n), c(-1) \big| u\in\U, c\in\CC, n\in\Z \}$ ($\U$, $\CC$ are some index sets).
Define the corresponding fields in $\L [[z,z^{-1}]]$:
$$ u(z) = \sum_{n\in\Z} u(n) z^{-n-1}, \quad c(z) = c(-1) z^0, \quad 
u\in\U, c\in\CC .$$
Let $\F$ be a subspace in $\L [[z,z^{-1}]]$ spanned by all the fields
$u(z), c(z)$ and their derivatives of all orders.

{\bf Definition.}
{\it
 A Lie algebra $\L$ with the basis as above is called a vertex Lie algebra
if the following two conditions hold:

(VL1) for all $u_1, u_2 \in \U$,
$$ [u_1(z_1), u_2(z_2) ] = \sum\limits_{j=0}^n f_j(z_2)
\left[ z_1^{-1} \left( {\d \over \d z_2} \right)^j \delta \left(
{z_2 \over z_1} \right) \right], \eqno{(\vla)}$$
where $f_j(z) \in\F, n \geq 0$ and depend on $u_1, u_2$,

(VL2) for all $c\in\CC$, the elements $c(-1)$ are central in $\L$.
}

\


Let $\L^{(+)}$ be a subspace in $\L$ with the basis $\{ u(n) \big| u\in\U, n\geq 0 \}$
and let  $\L^{(-)}$ be a subspace with the basis 
$\{ u(n), c(-1) \big| u\in\U, c\in\CC, n<0 \}$. Then $\L = \L^{(+)} \oplus \L^{(-)}$ and
$\L^{(+)}, \L^{(-)}$ are in fact subalgebras in $\L$.

The universal enveloping vertex algebra $V_{\L}$ of a vertex Lie algebra $\L$ 
is defined as an induced module
$$V_{\L} = \Ind_{\L^{(+)}}^\L (\C \o) = U(\L^{(-)}) \ot \o,$$
where $\C \o$ is a trivial 1-dimensional $\L^{(+)}$ module.

{\bf Theorem \voa. ([DLM], Theorem 4.8)} 
{\it
Let $\L$ be a vertex Lie algebra. Then

(a) $V_{\L}$ has a structure of a vertex algebra with the vacuum vector $\o$,
infinitesimal translation $D$ being a natural extension of 
the derivation of $\L$ given by 
$D(u(n))$  $=$ $-n u(n-1)$, $D(c(-1)) = 0$, $u\in\U$, $c\in\CC$,
and the state-field correspondence map $Y$ defined by the formula:
$$Y \left( a_1(-1-n_1) \ldots a_{k-1}(-1-n_{k-1}) a_k(-1-n_k) \o, z
\right) $$
$$ = 
:\left( {1\over n_1 !} \left( {\d \over \d z} \right)^{n_1} a_1 (z)
\right) 
\ldots
:\left( {1\over n_{k-1} !} \left( {\d \over \d z} \right)^{n_{k-1}} 
a_{k-1} (z) \right) 
\left( {1\over n_{k} !} \left( {\d \over \d z} \right)^{n_k} 
a_k (z) \right): \ldots : \quad 
, \eqno{(\Y)}$$
where $a_j \in \U, n_j \geq 0$ or $a_j \in\CC, n_j =0$.

(b) Any bounded $\L$-module is a vertex algebra module for $V_{\L}$.

(c) For an arbitrary character $\gamma: \CC \rightarrow \C$, the factor module
$$ V_{\L} (\gamma) =  U(\L^{(-)}) \o /  U(\L^{(-)}) \big< (c(-1) - \gamma(c)) 
\o \big>_{c\in\CC}$$
is a quotient vertex algebra.

(d) Any bounded $\L$-module in which $c(-1)$ act as $\gamma(c) \Id$, 
for all $c\in\CC$, is a vertex algebra module for $V_{\L}(\gamma)$.
}

In the formula (\Y) above,  
the normal
ordering of two fields $: a(z) b(z) :$ is defined as 
$$ : a(z) b(z) : = \sum\limits_{n < 0} a_{(n)} z^{-n-1} b(z) + 
\sum\limits_{n \geq 0} b(z) a_{(n)} z^{-n-1} .$$ 
Note the following relation that we will implicitly use throughout the paper:
$(a(-1) \o )_{(n)} = a(n)$ for $a \in \U$, \ $n\in \Z$.

\

{\bf Theorem \inv.} 
{\it Any $D$-invariant $\L$-submodule $U$ in $V_{\L}$ is a vertex algebra
ideal in $V_{\L}$. Conversely, every vertex algebra ideal in $V_{\L}$ is a 
$D$-invariant $\L$-submodule.
}

{\it Proof.} Let us prove the first part of the Theorem. We need to show that
for any $u\in U, \, a\in V_{\L}, \, n\in\Z$, we have $a_{(n)} u \in U$ and $u_{(n)} a \in U$.
By (\skeww), it is enough to prove that $a_{(n)} u \in U$. It is sufficient to consider
$a$ of the form $a = a_1(-1-n_1) \ldots a_k (-1-n_k) \o$, where
$a_j \in \U \cup \CC, \, n_j \geq 0$. We use induction on $k$. For $k=0$ we get that
$a = \o$ and $\o_{(n)} u = \delta_{n,-1} u$. The inductive step follows from the Borcherds'
formula (\Borb) and $\L$-invariance of $U$. The second part of the Theorem 
follows immediately from the definition of $V_{\L}$.   

\

{\bf Corollary \Cinv.} 
{\it
If $U$ is a maximal $D$-invariant $\L$-submodule in $V_\L$
then the quotient $L_\L = V_\L / U$ is a simple vertex algebra.
}

\

{\bf Remark \RVL.} In case when the set $\U$ contains an element $\omega$
generating the Virasoro field $\omega(z)$ in $\L$, satisfying
$[\omega(0), a(n)] = -n a(n-1)$ for all $a\in\U$,
the vertex algebra $V_{\L}$ becomes a VOA,
and in the statement of Theorem {\inv} the condition of $D$-invariance of $U$
will automatically follow from its $\L$-invariance.

\

{\bf 3.4. VOA associated with the twisted Virasoro-affine algebra.}

The toroidal VOA that will be constructed in Section 4, decomposes
into a tensor product of two VOAs. One of these factors is a VOA associated
with a twisted Virasoro-affine Lie algebra, which we introduce here.

Let $\df$ be a finite-dimensional reductive Lie algebra. Consider a 
semi-direct product of the Lie algebra of vector fields on a circle 
with a loop algebra:
$$ \tf = \Der \C[t_0, t_0^{-1}] \ltimes  \C[t_0, t_0^{-1}] \otimes \df .$$
A {\it twisted Virasoro-affine algebra} $\f$ is the universal central extension  
of the Lie algebra $\tf$.
Using the results on the central extensions of the loop algebras and 
the Lie algebra of vector fields on a circle, one can show that the second 
cohomology of $\tf$ has the following description:
$$ H^2(\tf) = S^2(\df)^\Inv \oplus \df^\Inv \oplus \C,$$
where the last 1-dimensional component corresponds to the Virasoro cocycle
on $\Der \C[t_0, t_0^{-1}]$, and its generator will be denoted by $C_\Vir$.
Since $\df$ is reductive, we have the following canonical projections of $\df$-modules:
$$\varphi: \quad \df \otimes \df \rightarrow S^2(\df)^\Inv$$
and
$$\psi: \quad \df \rightarrow Z(\df) = \df^\Inv .$$
We will use this maps to write down the Lie bracket in the twisted Virasoro-affine
algebra  $\f = \tf \oplus S^2(\df)^\Inv \oplus \df^\Inv \oplus \C$:
$$[L(n), L(m)] = (n-m) L(n+m) + {n^3 - n \over 12} \delta_{n,-m}  C_\Vir, \eqno{(\LL)}$$
$$[L(n), f(m)] = -m f(n+m) - (n^2 + n) \delta_{n,-m}  \psi(f), \eqno{(\Lf)}$$ 
$$[f(n), g(m)] = [f,g](n+m) + n \delta_{n,-m} \varphi(f \otimes g), \quad f,g \in \df. \eqno{(\fg)}$$ 
Here and below we are using the notations 
$L(n) = - t_0^{n+1} {d \over dt_0}$ and $f(n) = t_0^n \otimes f$ for $f \in \df$.

Consider the following fields in $\f$:
$$ \omega(z) = \sum_{n\in\Z} \omega(n) z^{-n-1} = \sum_{n\in\Z} L(n) z^{-n-2}$$
and
$$ f(z) = \sum_{n\in\Z} f(n) z^{-n-1}, \quad {\rm for} \quad f\in\df .$$

{\bf Proposition \tav.} 
{\it 
Twisted Virasoro-affine Lie algebra $\f$ is a vertex Lie algebra.
}

{\it Proof.} Take for a set $\U$ the element $\omega$ together with a basis of $\df$,
and for a set $\CC$ a basis of $S^2(\df)^\Inv \oplus \df^\Inv \oplus \C$.
Then the defining relations (\LL)-(\fg) may be rewritten as follows:
$$ \eqalign{
[\om(z_1), \om(z_2)] = \left( {\d \over \d z_2} \om(z_2) \right) \z 
&+ 2 \om(z_2) \zd \cr
&+ {C_\Vir \over 12} \zddd ,} \eqno{(\LLz)}$$
$$ \eqalign{
[\om(z_1), f(z_2)] = \left( {\d \over \d z_2} f(z_2) \right) \z 
&+  f(z_2) \zd \cr
&-  \psi(f)  \zdd ,} \eqno{(\Lfz)}$$
$$ [f(z_1), g(z_2)] = [f,g](z_2) \z +  \varphi(f \otimes g)  \zd . \eqno{(\fgz)}$$
This shows that $\f$ is indeed a vertex Lie algebra, and the claim of the Proposition is 
established.

 From now on we fix $\df$ to be $\df = \dg \oplus \glN$. 
 
 A linear map $S^2(\df)^\Inv \rightarrow \C$ defines a symmetric invariant bilinear form on $\df$.
 Thus $S^2(\df)^\Inv$ is the dual space to the space of symmetric invariant forms on $\df$ and has
 dimension 3. Let us fix an invariant form on $\dg$ normalized by the condition that 
 $(\alpha | \alpha) = 2$ for the long roots of $\dg$, an invariant form on $\slN$ with the same
 normalization, and a form on the space of scalar matrices normalized by its value on the identity
 matrix: $(I | I) = 1$. Denote by $\{ C_\sdg, C_\sln, C_\Hei \}$ the dual basis in $S^2(\df)^\Inv$.

The space $\df^\Inv \cong Z(\df)$ is one-dimensional, and we will denote its generator $\psi(I)$
by $C_\HV$. Hence the space $H^2(\tf)$ is 5-dimensional with the basis
$\{ C_\sdg, C_\sln, C_\Hei, C_\HV, C_\Vir \}$.

The Lie algebra $\f$ contains four subalgebras -- a Virasoro algebra, two affine algebras,
$\wdg = \C [t_0, t_0^{-1}] \otimes \dg \oplus \C C_\sdg$ and 
$\wsl = \C [t_0, t_0^{-1}] \otimes \slN \oplus \C C_\sln$,
and an infinite-dimensional Heisenberg algebra 
$\Hei = \C [t_0, t_0^{-1}] \otimes I \oplus \C C_\Hei$.

Fix a central character $\gamma: \quad H^2(\tf) \rightarrow \C$:
$$\gamma(C_\sdg) = c_\sdg, \quad \gamma(C_\sln) = c_\sln, \quad \gamma(C_\Hei) = c_\Hei, 
\quad \gamma(C_\HV) = c_\HV, \quad \gamma(C_\Vir) = c_\Vir ,$$
and consider the corresponding quotient $V_\sf (\gamma)$ of 
the universal enveloping vertex algebra.
Using the commutator formula (\comm), we derive from (\Lfz), (\fgz) 
the following relations for the $n$-th products.

{\bf Lemma \reln.} 
{\it
The following relations hold in $V_\sf (\gamma)$:

$$\eqalign{
(a) \quad\quad\quad E_{ab} (0) E_{cd} (-1) \o &= \delta_{bc} E_{ad} (-1) \o - \delta_{ad} E_{cb} (-1) \o, 
 \cr
E_{ab} (1) E_{cd} (-1) \o &= \delta_{ad} \delta_{bc} c_\sln \o 
+ \delta_{ab} \delta_{cd} \left( {c_\Hei \over N^2} - {c_\sln \over N} \right) \o, \cr
 E_{ab} (n) E_{cd} (-1) \o &= 0 \hbox{\rm \ for \ } n\geq 2. \cr}$$

$$\eqalign{
(b) \quad\quad \om_{(0)} E_{ab} (-1) \o &= D (E_{ab} (-1) \o), \quad
\om_{(1)} E_{ab} (-1) \o = E_{ab} (-1) \o, \cr
\om_{(2)} E_{ab} (-1) \o &= - \delta_{ab} {2 c_\HV \over N} \o, \quad
\om_{(n)} E_{ab} (-1) \o = 0 \hbox{\rm \ for \ } n\geq 3.\cr}$$
}





\

Let us now discuss bounded weight modules for $\f$. This Lie algebra is $\Z$-graded
by degree in $t_0$. We associate with this grading a decomposition $\f = \f_- \oplus
\f_0 \oplus \f_+$, where $\f_0 = \C d_0 \oplus \dg \oplus \glN \oplus H^2(\tf)$.
Let $V$ be a finite-dimensional irreducible $\dg$-module, and $W$ be a finite-dimensional  
irreducible module for $\slN$. Fix a central character $\gamma: H^2(\tf) \rightarrow \C$
and two constants $h_\Vir, h_\Hei \in \C$. We define on $V\otimes W$ the structure
of an irreducible $\f_0$-module on which $L(0) = -d_0$ acts as multiplication by $h_\Vir$,
$I$ acts as multiplication by $h_\Hei$, and the action of $H^2(\tf)$ is determined by $\gamma$.
Let $\f_+$ act on $V\otimes W$ trivially and consider the induced module
$$M_\sf (V,W,h_\Hei,h_\Vir,\gamma) = \Ind_{\sf_0\oplus \sf_+}^\sf \left(V \otimes W \right) .$$
This module has a unique maximal submodule, and the factor-module by the maximal submodule
is an irreducible $\f$-module, which we denote as $L_\sf (V,W,h_\Hei,h_\Vir,\gamma)$.

{\bf Remark \RVA.} Similar constructions may be applied to Virasoro, affine, and Heisenberg
algebras, yielding the corresponding vertex algebras
 $V_\Vir (c_\Vir)$, $V_\swdg (c_\sdg)$, $V_\wsl (c_\sln)$, $V_\Hei (c_\Hei)$ and
irreducible highest weight modules $L_\Vir (h_\Vir, c_\Vir)$ for the Virasoro
algebra, $L_\swdg (V, c_\sdg)$ for the affine algebra $\wdg$, $L_\wsl (W,c_\sln)$ for the affine algebra
$\wsl$ and $L_\Hei (h_\Hei, c_\Hei)$ for the infinite-dimensional Heisenberg algebra.

Note that for the trivial 1-dimensional modules $V = \C$, $W = \C$ and $h_\Vir = h_\Hei = 0$,
the irreducible module $L_\sf (\C, \C, 0, 0, \gamma)$ is precisely the simple VOA
$L_\sf (\gamma)$.

For a generic $\gamma$ ($\gamma$ not at a critical level), we may apply the Sugawara 
construction to decompose the irreducible module $L_\sf (V,W,h_\Hei,h_\Vir,\gamma)$
into a tensor product of irreducible Virasoro, affine and Heisenberg modules.

{\bf Proposition \FA.}
Let $c_\sdg \neq - h^{\vee}, c_\sln \neq -N, c_\Hei \neq 0$, where $h^{\vee}$ is 
the dual Coxeter number for $\wdg$.
Then the VOA $V_\sf (\gamma)$ decomposes into a tensor product of four VOAs:
$$V_\sf (\gamma) \cong V_\swdg (c_\sdg) \otimes V_\wsl (c_\sln) \otimes V_\Hei (c_\Hei) 
\otimes V_\Vir (c^\prime_\Vir) ,$$
where 
$$c^\prime_\Vir = c_\Vir - {c_\sdg \dim (\dg) \over c_\sdg
+ h^\vee} - {c_\sln (N^2 -1) \over c_\sln + N}
-1 + 12 {c_\HV^2 \over c_\Hei}, \eqno{(\cvp)}$$
and the Heisenberg VOA $V_\Hei (c_\Hei)$ is taken with a non-standard Virasoro element
$$\om_\Hei = {1 \over 2 c_\Hei} I(-1) I(-1) \o + {c_\HV \over c_\Hei} I(-2) \o, \eqno{(\vhei)}$$
so that its rank is $1 - 12 {c_\HV^2 \over c_\Hei}$.

{\it Proof.} This result is obtained by applying the Sugawara construction three times -- to the affine 
$\wdg$-subalgebra,
affine $\wsl$-subalgebra, and the twisted Virasoro-Heisenberg subalgebra (see e.g. [FLM] and [ACKP] for details).
Since this construction is well-known, we only sketch the proof.

Let $\{ u_i \}, \{ u^i \}$ be dual bases of $\dg$, and  $\{ v_j \}, \{ v^j \}$ be dual bases of $\sln$ with 
respect to the chosen invariant bilinear forms.
 Consider a new Virasoro field
$$\eqalign{
\omega^\prime(z) =& \omega(z) - {1 \over 2(c_\sdg + h^{\vee})} \sum_i : u_i(z) u^i (z) : 
- {1 \over 2(c_\sln + N)} \sum_j : v_j(z) v^j (z) : \cr
& - {1\over 2 c_\Hei} : I(z) I(z) : 
- {c_\HV \over c_\Hei} {\d \over \d z} I(z)  .} \eqno{(\newVir)}$$
It is possible to verify that the moments of $\omega^\prime(z)$ satisfy the Virasoro algebra relations, 
with the action of the central element given by (\cvp). 
Moreover, this new Virasoro field $\omega^\prime (z)$ commutes with the fields of the affine
$\wdg$ and $\wsl$ subalgebras, as well the Heisenberg subalgebra field. 

The formula (\newVir) defines a homomorphism of vertex algebras
$$  V_\swdg (c_\sdg) \otimes V_\wsl (c_\sln) \otimes V_\Hei (c_\Hei) 
\otimes V_\Vir (c^\prime_\Vir) \rightarrow V_\sf (\gamma) ,$$
which is in fact an isomorphism. Moreover, if we choose the Virasoro element in $V_\Hei (c_\Hei)$
to be given by (\vhei), the above map becomes the isomorphism of the VOAs.

\

{\bf Corollary \FAA.}
{\it
Under the same restriction on the central charges as in Proposition \FA, 
the irreducible highest weight $\f$-module
$L_\sf (V,W,h_\Hei, h_\Vir, \gamma)$ decomposes into a tensor product of irreducible highest weight
modules for the affine $\wdg$, $\wsl$,
infinite-dimensional Heisenberg and the Virasoro modules:
$$L_\sf (V,W,h_\Hei, h_\Vir, \gamma) \cong 
L_{\swdg} (V, c_\sdg) \otimes L_{\wsl} (W, c_\sln) \otimes L_{\Hei} (h_\Hei, c_\Hei) 
\otimes L_{\Vir} (h_\Vir^\prime, c_\Vir^\prime),$$
where $c_\Vir^\prime$ is given by (\cvp) and
$$h_\Vir^\prime = h_\Vir - {\Omega_V \over 2(c_\sdg + h^\vee)} - {\Omega_W \over 2(c_\sln + N)}
 - {h_\Hei^2 - 2c_\HV h_\Hei \over 2 c_\Hei} . \eqno{(\hvp)}$$
}

Here $\Omega_V$ and $\Omega_W$ are the eigenvalues of the Casimir operators of $\dg$ and $\slN$ on $V$
and $W$ respectively. These are given by $\Omega_V = (\lambda_V | \lambda_V + 2\rho)$,
$\Omega_W = (\lambda_W | \lambda_W + 2\rho)$, where $\lambda_V$ is the highest weight of the irreducible
$\dg$-module $V$, and $\lambda_W$ is the highest weight of the $\slN$-module $W$ [K1].

\

{\bf Remark \Rp.} If one of the inequalities in the statement of the above proposition fails, 
we still  can apply a partial Sugawara construction to the remaining components. 
For example, if  $c_\sdg \neq - h^{\vee}, c_\sln \neq -N$, but $c_\Hei =0$,
the irreducible highest weight $\f$-module is isomorphic to the tensor product of the affine 
$\wdg$ and  $\wsl$-modules and an irreducible highest weight module for the twisted 
Heisenberg-Virasoro algebra at level zero.
The characters of such modules for the twisted Heisenberg-Virasoro algebra were computed in [B3].


\

As we mentioned earlier, the toroidal VOA decomposes in a tensor product, where one
factor is a twisted Virasoro-affine VOA. The other factor in this decomposition is
a sub-VOA of a lattice VOA, which we will describe next.

\

{\bf 3.5. Hyperbolic lattice VOA.}

Here we present the construction of a hyperbolic lattice VOA. The general 
construction of a VOA corresponding to an arbitrary even lattice may be found 
in [FLM] or [K2]. 

Consider a hyperbolic lattice $\Hyp$, which is a free abelian group on $2N$
generators 
\break
$\{ u_i , v_i |  i = 1, \ldots, N \}$ with the symmetric bilinear
form
$$ ( \cdot | \cdot ) : \quad \Hyp \times \Hyp \rightarrow \Z ,$$
defined by 
$$ (u_i | v_j) = \delta_{ij} , \quad (u_i | u_j) = (v_i | v_j) = 0.$$
Note that the form $(\cdot | \cdot)$ is non-degenerate and $\Hyp$ is an 
even lattice, i.e., $(x | x) \in 2\Z$.

The construction of the VOA associated to $\Hyp$ proceeds as follows.

 First we complexify $\Hyp$:
$$ H = \Hyp \ot_{\Z} \C,$$
and extend $(\cdot | \cdot)$ by linearity on $H$. Next, we ``affinize'' $H$
by defining a Lie algebra
$\widehat H = \C[t, t^{-1}] \ot H \oplus \C K$
with the bracket
$$ [x(n), y(m)] = n (x | y) \delta_{n, -m} K, \quad x,y \in H, 
\quad [\widehat H, K] = 0. \eqno{(\xyK)} $$
Here and in what follows, we are using the notation $x(n) = t^n \ot x$. The algebra
$\widehat H$ has a triangular decomposition $\widehat H = \widehat H_-
\oplus \widehat H_0 \oplus \widehat H_+$, where 
$\widehat H_0 = \left< 1\ot H, K \right>$ and $\widehat H_\pm = 
t^{\pm 1} \C [t^{\pm 1}] \ot H$.

We also need a twisted group algebra of $\Hyp$, denoted by $\C[\Hyp]$, which
we now describe. The basis of $\C[\Hyp]$ is $\{ e^x | x \in \Hyp \}$, 
and the multiplication is twisted with the 2-cocycle $\epsilon$:
$$ e^x e^y = \epsilon(x,y) e^{x+y} , \quad x,y \in \Hyp , \eqno{(\twi)}$$
where $\epsilon$ is a multiplicatively bilinear map
$$\epsilon: \Hyp \times \Hyp \rightarrow \{ \pm 1 \},$$
defined on the generators by $\epsilon(v_i, u_j) = (-1)^{\delta_{ij}},
\epsilon(u_i, v_j) = \epsilon(u_i, u_j) = \epsilon(v_i, v_j) = 1,
\quad i,j = 1,\ldots, N$.

 We define the structure of $\widehat H_0 \oplus \widehat H_+$-module on 
$\C[\Hyp]$, letting $\widehat H_+$ act on $\C[\Hyp]$ trivially and
$\widehat H_0$ act by
$$ x(0) e^y = (x | y) e^y, \quad K e^y = e^y. \eqno{(\xoe)}$$

 Finally let $V_\hyp$ be the induced $\widehat H$ module:
$$ V_\hyp = \Ind_{\widehat H_0 \oplus \widehat H_+}^{\widehat H} \left(
\C[\Hyp] \right) .$$
This is the VOA attached to the lattice $\Hyp$. As a space $V_\hyp$ is 
isomorphic to the tensor product of the symmetric algebra $S(\widehat H_-)$
with the twisted group algebra $\C[\Hyp]$:
$$ V_\hyp = S(\widehat H_-) \ot \C[\Hyp] .$$

 The $Y$-map is defined on the basis elements of $\C[\Hyp]$ by
$$ Y(e^x, z) = \exp \left( \sum\limits_{j \geq 1} {x(-j) \over j} z^j \right)
  \exp \left( - \sum\limits_{j \geq 1} {x(j) \over j} z^{-j} \right)
e^x z^x , \eqno{(\Ye)}$$
where $e^x$ acts by twisted multiplication (\twi) and $z^x e^y = z^{(x|y)} e^y$.
For a general basis element $a = x_1 (-1-n_1) \ldots x_k(-1-n_k) \ot e^y$,
with $x_i, y \in \Hyp, n_i \geq 0$, one defines (cf. (\Y))
$$Y(a,z) = :\left( {1\over n_1 !} \left({\d \over \d z}
\right)^{n_1} x_1(z) \right) \ldots
\left( {1\over n_k !} \left({\d \over \d z}
\right)^{n_k} x_k(z) \right) Y(e^y,z) :, \eqno{(\Yw)}$$
where $x(z) = \sum\limits_{j\in\Z} x(j) z^{-j-1}$.

The Virasoro element in $V_\hyp$ is $\omega_\hyp = \sum\limits_{p=1}^N
u_p (-1) v_p(-1) \o$, where $\o = e^0$ is the identity element
of $V_\hyp$. The rank of $V_\hyp$ is $2N$.

In the construction of the toroidal VOAs we would need not $V_\hyp$ itself,
but its sub-VOA $V_\hyp^+$:
$$  V_\hyp^+ = S(\widehat H_-) \ot \C[\Hyp^+] ,$$
where $\Hyp^+$ (resp. $\Hyp^-$) is the isotropic sublattice of $\Hyp$ 
generated by $\{ u_i |  i = 1, \ldots, N \}$ (resp. 
$\{ v_i |  i = 1, \ldots, N \}$). 
 One can verify immediately by inspecting
(\Ye) and (\Yw) that $V_\hyp^+$ is indeed a sub-VOA of $V_\hyp$. 
Also note that
the cocycle $\epsilon$ trivializes on $\C[\Hyp^+]$, making $\C[\Hyp^+]$
the usual (untwisted) group algebra.
The Virasoro element of $V_\hyp^+$ is the same as in $V_\hyp$, and so
the rank of $V_\hyp^+$ is also $2N$.

Let us describe a class of modules for $V_\hyp^+$. Consider the group
algebra $\C[H^+]$ of the vector space $H^+ = \Hyp^+ \ot_{\Z} \C.$
The space $S(\widehat H_-) \ot \C[H^+] \ot \C[\Hyp^-]$ has a 
structure of a VOA module for $V_\hyp^+$, where the action of 
$V_\hyp^+$ is still given by (\xoe),(\Ye) and (\Yw). 
Fix $\alpha \in \C^N, \beta\in\Z^N$. Then the subspace
$$ M_\hyp^+(\alpha,\beta) = S(\widehat H_-) \ot 
e^{\alpha \u + \beta \vv } \C[\Hyp^+] $$
in $S(\widehat H_-) \ot \C[H^+] \ot \C[\Hyp^-]$ is an irreducible
VOA module for $V_\hyp^+$. Here we are using the notations
$\alpha \u = \alpha_1 u_1 + \ldots +  \alpha_N u_N$, etc.
For $\beta = 0$ we will denote the module $M_\hyp^+(\alpha,0)$ simply
by $M_\hyp^+(\alpha)$.

\

\

{\bf 4. Toroidal vertex operator algebras.}

\

In this section we will construct several VOAs associated with the toroidal Lie algebras.
We will construct 
a universal enveloping VOA $V_\sg$, its ``level $c$'' quotient $V_\sg (c)$ and
the simple quotient $L(T_0)$. As $\g$-modules,
$V_\sg (c)$ does not belong to category $\B_\chi$, but $L(T_0)$ does. We will establish several 
important relations that hold in $L(T_0)$, which will allow us to show that
$L(T_0)$ factors into the tensor
product of two VOAs discussed in the previous section, $V_\Hyp^+$ and the twisted
Virasoro-affine VOA $L_\sf (\gamma_0)$.


 The key fact which makes it possible to construct these VOAs, is the 
observation that toroidal Lie algebras $\g(\mu,\nu)$ are in fact vertex Lie algebras for
all values of $\mu, \nu$. 

 This observation is not quite trivial, since it requires a rather delicate choice of a
basis in $\g(\mu,\nu)$, in order to exhibit the vertex Lie algebra structure.

\

{\bf Theorem {\tvla}.} 
{\it Toroidal Lie algebras $\g(\mu,\nu)$ are vertex Lie algebras.}

{\it Proof.}
 Consider the following generating series in $\g [[z,z^{-1}]]$:
$$ k_0 (\r,z ) = \sum_{j=-\infty}^{\infty} t_0^j t^\r k_0 z^{-j} , \quad\quad  
 k_p (\r,z ) = \sum_{j=-\infty}^{\infty} t_0^j t^\r k_p z^{-j-1} ,\eqno{(\korz)}$$ 
$$ g(\r,z ) = \sum_{j=-\infty}^{\infty} t_0^j t^\r g z^{-j-1} , \quad g\in\dg, \eqno{(\grz)}$$ 
$$ \td_p (\r,z ) = \sum_{j=-\infty}^{\infty} t_0^j t^\r \td_p z^{-j-1} ,\quad\quad
 \td_0 (\r,z ) = \sum_{j=-\infty}^{\infty} t_0^j t^\r \td_0 z^{-j-2} ,\eqno{(\dorz)}$$ 
where for $p=1,\ldots,n,$
$$ t_0^j t^\r \td_p = t_0^j t^\r d_p - \nu r_p t_0^j t^\r k_0 , \eqno{(\tdp)}$$
and
$$ t_0^j t^\r \td_0 = - t_0^j t^\r d_0 + (\mu + \nu) (j+{1\over 2}) t_0^j t^\r k_0 . \eqno{(\tdo)}$$

Although the moments of the above series are not linearly independent, all linear
dependencies may be encoded as relations between the fields:
$$ {\d \over \d z} k_0 (\r,z) = \sum\limits_{p=1}^N r_p k_p (\r, z). $$
Using these relations we can eliminate from the above list for each non-zero $\r$ the 
field $k_p(\r,z)$ with the smallest $p$ such that $r_p \neq 0$. Non-zero moments of the 
remaining fields will form a basis of $\g(\mu,\nu)$.

Verification of the axioms of a vertex Lie algebra is now quite straightforward. 
The set $\CC$ consists of a single element that corresponds to the central field 
$k_0(0,z) = k_0 z^0$, so (VL2) holds.

Before we check the property (VL1), let us record the commutator relations between 
the newly introduced elements $t_0^j t^\r \td_p, t_0^j t^\r \td_0$. Note that their
brackets with the elements of $\R\otimes \dg$ and $\K$ are essentially given by (\Ldg) and (\Ldk)
(with an obvious change of sign for $t_0^j t^\r \td_0$),
while the rest of the commutator relations are given by the following formulas with $a,b = 1,\ldots, N$:
$$ \eqalign{ [t_0^i t^r \td_a, t_0^j t^s \td_b] = 
& s_a t_0^{i+j} t^{r+s} \td_b - r_b t_0^{i+j} t^{r+s} \td_a \cr 
+ & (\mu s_a r_b + \nu r_a s_b) j t_0^{i+j} t^{r+s} k_0
+ (\mu s_a r_b + \nu r_a s_b) \sum_{p=1}^N s_p t_0^{i+j} t^{r+s} k_p . \cr} \eqno{(\tdadb)}$$
$$ \eqalign{ [t_0^i t^r \td_0, t_0^j t^s \td_b] = & - j t_0^{i+j} t^{r+s} \td_b - r_b t_0^{i+j} t^{r+s} \td_0
 - (\mu r_b (j-1)  + \nu s_b (i+1) )  j t_0^{i+j} t^{r+s} k_0 \cr
- & (\mu r_b j  + \nu s_b (i+1) ) \sum_{p=1}^N s_p t_0^{i+j} t^{r+s} k_p . \cr} \eqno{(\tdodb)}$$
$$ \eqalign{ [t_0^i t^r \td_0, t_0^j t^s \td_0] = & (i-j) t_0^{i+j} t^{r+s} \td_0
 + (\mu+\nu) j(j+1)(i+1) t_0^{i+j} t^{r+s} k_0 \cr
+ & (\mu+\nu) (j+1)(i+1) \sum_{p=1}^N s_p t_0^{i+j} t^{r+s} k_p . \cr} \eqno{(\tdodo)}$$

Using these formulas together with (\Ldg) and (\Ldk), 
we can derive the commutator relations for the fields in $\g$:

$$ \left[ k_a (\r, z_1), k_b (\m, z_2) \right] = 0, \eqno{(\Rkakb)}$$
$$ \left[  g(\r, z_1), k_a (\m, z_2) \right] = 0, \eqno{(\Rgka)}$$
$$ \left[  g_1(\r, z_1), g_2 (\m, z_2) \right] = [g_1,g_2](z_2) \z  $$
$$  + (g_1 | g_2) k_0(\r+\m, z_2) \zd 
+ (g_1 | g_2) \sum_{p=1}^N r_p k_p (\r+\m, z_2) \z ,\eqno{(\Rgg)}$$
$$ \left[ \td_j (\r, z_1), g (\m, z_2) \right] = m_j g(\r+\m,z_2) \z ,\eqno{(\Rdjg)}$$
$$ \left[ \td_0 (\r, z_1), g (\m, z_2) \right] = \dzb \left\{ g(\r+\m, z_2) \z \right\} ,\eqno{(\Rdog)}$$ 
$$\eqalign{ \left[ \td_i (\r, z_1), \td_j (\m, z_2) \right] = &
\left( m_i \td_j (\r+\m,z_2) - r_j \td_i (\r+\m,z_2) \right) \z \cr
 - & (\mu m_i r_j + \nu r_i m_j) \sum_{p=1}^N r_p k_p (\r+\m, z_2) \z \cr
 - & (\mu m_i r_j + \nu r_i m_j) k_0 (\r+\m,z_2) \zd , \cr} \eqno{(\Rdidj)}$$
$$ \left[ \td_0 (\r, z_1), \td_j (\m, z_2) \right] = \dzb \left\{ \td_j(\r+\m,z_2) \z \right\}
%
 - r_j \td_0 (\r+\m, z_2) \z $$
$$\quad\quad\quad + \nu m_j \sum_{p=1}^N r_p k_p (\r+\m,z_2) \zd 
 + \nu m_j k_0 (\r+\m,z_2) \zdd $$
$$ - \mu r_j \dzb \left\{ \sum_{p=1}^N r_p k_p (\r+\m,z_2) \z 
%
+ k_0 (\r+\m,z_2) \zd \right\} ,\eqno{(\Rdodj)}$$
$$ \left[ \td_0 (\r, z_1), \td_0 (\m, z_2) \right] = \left\{ \dzb \td_0 (\r+\m,z_2) \right\} \z
+ 2 \td_0 (\r+\m,z_2) \zd $$
$$ + (\mu+\nu) \dzb \left\{ \sum_{p=1}^N r_p k_p (\r+\m,z_2) \zd 
+ k_0 (\r+\m,z_2) \zdd \right\} ,\eqno{(\Rdodo)}$$
where \quad $g,g_1, g_2 \in \dg, \quad a,b = 0,1,\ldots,N , \quad i,j=1,\ldots,N$.

 Now, the right-hand sides of the above commutators are precisely in the format required by (VL1).
 Thus both (VL1) and (VL2) hold, we conclude that $\g(\mu,\nu)$ is a vertex Lie algebra,
 which allows us to consider its universal enveloping vertex algebra $V_\sg$.
Moreover, $\td_0 (0,z)$ is the Virasoro field with the central element
$C_\Vir = 12 (\mu + \nu) k_0$ and $D = t_0^{-1} \td_0.$ 

 The subalgebra $\g^{(-)}$ of the vertex Lie algebra $\g$ (not to be confused
with its subalgebra $\g_-$ with respect to $\Z$-grading) is generated by the following elements: 
 $t_0^j t^r k_0$ with $j\leq 0$, $t_0^j t^r k_p, \, t_0^j t^r g, \, t_0^j t^r \td_p$
with $j\leq -1, p=1,\ldots,N,$ and $t_0^j t^r \td_0$ with $j\leq -2$. The subalgebra 
$\g^{(+)}$ is spanned by the complementary moments of the fields (\korz)-(\dorz).
We recall that $\g^{(+)} \o = 0$ in $V_\sg$.

 It follows from Theorem {\voa} that $Y(t^r k_0, z) = k_0 (r,z), \;$ 
 $Y(t_0^{-1}t^r k_j, z) = k_j (r,z), $
 \break 
 $Y(t_0^{-1}t^r g, z) = g (r,z), \;$ 
 $Y(t_0^{-1}t^r \td_j, z) = \td_j (r,z), \;$
 $Y(t_0^{-2}t^r \td_0, z) = \td_0 (r,z)$.
For the sake of simplicity of notations, we are writing $Y(t^r k_0,z)$ for $Y( (t^r k_0) \o, z)$, etc. 
Also, when $r=0$, we will simply write $g(z)$ for $g(0,z)$, etc.    

We are going to show that for a particular irreducible $\g_0$-module $T_0$, the irreducible $\g$-module 
$L(T_0)$ is a factor-VOA of $V_\sg$. This will be done in two steps. First we will
construct a factor-VOA of $V_\sg$ that has an irreducible $\g_0$-module as its top. After that
we will show that the irreducible quotient of this $\g$-module is a VOA, and determine the structure 
of this vertex algebra.

The operator $\td_0 =
- d_0 + {1\over 2} (\mu+\nu) k_0$ induces the $\Z$-grading of the universal enveloping 
vertex algebra $V_{\L}$.
Since $V_\sg = U(\g^{(-)}) \otimes \o$, we see that its zero component is spanned by the elements
$(t^{r_1} k_0) \ldots (t^{r_s} k_0) \o$.  

{\bf Proposition \tops.}
Fix a non-zero $c\in\C$ and consider a $\g$-submodule $R(S)$ in $V_\sg$ generated by the set
$ S = \left\{ k_0 \o - c\o, \; (t^r k_0)(t^m k_0) \o - c (t^{r+m} k_0) \o \; | \; r,m \in \Z^N \right\} $. 

(a) The quotient $V_\sg (c) = V_\sg / R(S)$ is a factor-module of the generalized Verma module
$M(T_0)$ with the top
$$ T_0 = \C [q_1^\pm, \ldots, q_N^\pm] \otimes V \otimes W = \C [q_1^\pm, \ldots, q_N^\pm]$$
defined as in (\TT) with $\alpha = 0$, $d = {1\over 2} (\mu + \nu) c$,
 $V$ being the trivial 1-dimensional $\dg$-module,
$W$ being a 1-dimensional $\glN$-module on which $\slN$ acts trivially and $I$ acts a multiplication by 
$h = N \nu c$. 

(b) $V_\sg (c)$ inherits a vertex algebra structure from $V_\sg$.

(c) $V_\sg (c)$ is a VOA of rank $12 (\mu + \nu) c$ with the Virasoro field $\om(z) = \td_0 (z)$.

(d) The projection from $M(T_0)$ to $L(T_0)$ factors through $V_\sg (c)$:
$$ M(T_0) \rightarrow V_\sg(c) \rightarrow L(T_0) . \eqno{(\diag)}$$
This defines a VOA structure on $L(T_0)$ as a factor-VOA of $V_\sg(c)$.

{\it Proof.} Let us prove part (a). Consider a $\Z$-grading on $V_\sg$. We claim that
the zero component $R(S)_0$ coincides with 
$$ \Span \left< (t^{r_1} k_0) \ldots (t^{r_s} k_0) \o - c^{s-1} (t^{r_1+\ldots + r_s} k_0) \o,
 \quad k_0 \o - c\o \; | \; r_1, \ldots, r_s \in\Z^N \right> .
\eqno{(\RS)}$$
First of all it is easy to see by repeated multiplication of the elements in $S$ by
$t^{r_j} k_0$, that the elements (\RS) are indeed in $R(S)$. Next, we write
$R(S) = U(\g) S$. We have a triangular decomposition $U(\g) = U(\g_-) \otimes U(\g_0)
\otimes U(\g_+)$, and $\g_+$ acts on $S$ trivially because the elements of $S$ are of degree zero. 
This implies that $R(S)_0 = U(\g_0) S$. Let us show that (\RS) is invariant
under the action of $\g_0$. The subalgebra $\g_0$ is spanned by the elements $t^m k_0, \; t^m k_p, \;
t^m g, \; t^m \td_p, \; t^m \td_0, \quad m\in\Z^N, g \in \dg, p=1,\ldots, N$. We have already verified
the invariance of (\RS) under the action of $t^m k_0$. We note that the remaining
generators of $\g_0$ belong to $\g^{(+)}$ and thus annihilate $\o$. In addition to this,
the elements $t^m k_p$ and $t^m g$ commute with $t^r k_0$, which implies that 
they also annihilate (\RS).
It follows from (\Ldk) that $t^m \td_0$ annihilate (\RS) as well.

 Using the commutator relation $\left[ t^m \td_p, t^r k_0 \right] = r_p t^{r+m} k_0$, we get that
$$ (t^m \td_p) \left( (t^{r_1} k_0) \ldots (t^{r_s} k_0) \o - c^{s-1} (t^{r_1+\ldots + r_s} k_0) \o 
\right) $$
$$= \sum_{j=1}^s r_j \left( (t^{r_1} k_0) \ldots (t^{r_j+m} k_0) \ldots  (t^{r_s} k_0) \o 
- c^{s-1} (t^{r_1+\ldots + r_s + m } k_0) \o \right) ,$$
and the right hand side is in (\RS). This proves our claim that $R(S)_0$ is given by (\RS).

It follows from this that the top of the module $V_\sg (c) = V_\sg / R(S)$ may be identified
with the space of Laurent polynomials $T_0 = \C [q_1^\pm,\ldots, q_N^\pm]$, under the isomorphism 
$$ (t^r k_0) \o \mapsto c q^r .$$
Let us describe the action of the subalgebra $\g_0$ on the top $T_0$. It follows from 
the relations (\RS) that
$$ (t^m k_0) q^r = c q^{r+m} . \eqno{(\koq)}$$
Next, we have seen that $t^m k_p, t^m g$ and $t^m \td_0$ annihilate $T_0$:
$$ (t^m k_p) q^r = 0, \quad (t^m g) q^r = 0, \quad (t^m \td_0) q^r = 0. \eqno{(\kpq)} $$
 Since $t^m \td_0 = - t^m d_0 + {1\over 2} (\mu + \nu) t^m k_0$, we get that
$$ (t^m d_0) q^r = {1\over 2} (\mu + \nu) c q^{r+m} . \eqno{(\doq)}$$
Finally,
$$ (t^m \td_p) q^r = {1\over c} (t^m \td_p)  (t^r k_0) \o =  
{1\over c} \left[t^m \td_p, t^r k_0 \right] \o =
{1\over c} r_p (t^{r+m} k_0) \o = r_p q^{r+m}. \eqno{(\tdpq)} $$
Taking into account that $t^m \td_p = t^m d_p - m_p \nu t^m k_0$, we obtain
$$ (t^m d_p) q^r = (r_p + \nu c m_p) q^{r+m} , \eqno{(\dpq)}$$
which corresponds to the tensor module action (\dvw) with the trivial action of $\slN$
and $I$ acting as multiplication by $N\nu c$.
This completes the proof of part (a). Part (b) follows from Theorem \inv.
The claim of part (c) has been already established and 
finally, for part (d) we note that $L(T_0)$ is a unique irreducible factor of $M(T_0)$, thus 
the projection $M(T_0) \rightarrow L(T_0)$ factors through $V_\sg(c)$ as in (\diag).
Let us point out that the kernel of the homomorphism $M(T_0) \rightarrow V_\sg (c)$
is the submodule generated by $\left\{ (t_0^{-1} t^m \td_0) \o |
m \in \Z^N \right\} $. Applying Theorem {\inv} again, we conclude that $L(T_0)$ inherits the vertex 
operator algebra structure from $V_\sg (c)$ given by formula (\Y).
 This completes the proof of the Proposition.

 Next we are going to study the structure of the VOA $L(T_0)$.

{\bf Theorem \TA.} 

{\it
(a) The VOA $L(T_0)$ is generated by the following elements: $q^m = {1\over c} (t^m k_0) \o$, $(t_0^{-1} g) \o$, 
$(t_0^{-1} k_a) \o$ , $(t_0^{-1} \td_a) \o$, $E_{ab}, (t_0^{-2} \td_0) \o$, with $m\in\Z^N, g\in\dg,
a,b = 1, \ldots, N$, where 
$$ E_{ab} = {1 \over c} (t_0^{-1} t^{\epsilon_a} \td_b) (t^{-\epsilon_a} k_0) \o
- (t_0^{-1} \td_b) \o +   {1 \over c} \delta_{ab} (t_0^{-1} k_a) \o, \eqno{(\Eab)}$$
with $\epsilon_a$ being a standard basis vector in $\Z^N$ with $1$ in $a$-th position.

(b) The module $L(T_0)$ has a structure of a $\C[q_1^\pm, \ldots, q_N^\pm]$-module, which is
compatible with the action of the algebra of Laurent polynomials on $T_0$. 
The action of the field $k_0 (m,z)$ is given by the following vertex operator:
$${1\over c} k_0 (m,z)  =  Y(q^m,z) =  
q^m \exp\left( \sum_{p=1}^N m_p \sum_{j\in\Z \backslash \{ 0 \} } {1\over cj} (t_0^{-j} k_p) z^j 
\right) . \eqno{(\Yko)}$$

(c) The action on $L(T_0)$ of the remaining fields in
(\korz)-(\dorz) is expressed in the following way
$$  g(m,z) =  g(z)  Y(q^m, z), \eqno{(\Ygm)}$$
$$ k_a(m , z) = k_a(z)  Y( q^m, z), \eqno{(\Ykam)}$$
$$  \td_a(m , z) = :\td_a(z)  Y(q^m, z): 
+ \sum_{p=1}^N m_p Y(E_{pa},z) Y( q^m, z)  , \eqno{(\Ydam)}$$
$$ \td_0(m, z) =  :\td_0(z)  Y(q^m, z): 
+ {1\over c} \sum_{a,b=1}^N m_a k_b(z) Y(E_{ab},z) Y( q^m, z)$$
$$+ (\mu - {1\over c}) \sum_{p=1}^N m_p \left( \dz k_p (z) \right)  Y( q^m, z) . \eqno{(\Ydom)}$$


(d) The vertex operator algebra $L(T_0)$ is isomorphic to the
tensor product of two VOAs:
$$ L(T_0) = \VH \otimes L_\sf (\gamma_0) ,$$
where $\VH$ is a sub-VOA of the lattice VOA described in section 3.5, 
and $L_\sf(\gamma_0)$ is the simple VOA corresponding to
the twisted Virasoro-affine Lie algebra constructed from the reductive Lie algebra $\df = \dg \oplus \glN$,
and the central character $\gamma_0$ given by the following values:
$$ c_\sdg = c, \quad \quad c_\sln = 1 - \mu c, $$
$$ c_\Hei = N(1-\mu c) - N^2 \nu c, \quad\quad c_{\HV} = N({1 \over 2} - \nu c), $$
$$c_\Vir = 12 c (\mu + \nu) - 2N. \eqno{(\gamo)}$$ 
The fields $k_0 (m,z)$, $k_p(z)$, $\td_p (z)$, $\;  p = 1, \ldots, N$, act on $\VH$ by
$$k_0 (m,z) = c Y(e^{mu},z), \quad k_p(z) = c u_p (z), \quad \td_p(z) = v_p (z) ,$$
while the fields $g(z)$ and $Y(E_{ab},z)$ act on $L_\sf (\gamma_0)$. The field $\td_0 (z)$
is the Virasoro field of the tensor product $\VH \otimes L_\sf (\gamma_0)$.
}


\


{\it Proof.}
 We are going to determine the structure of the VOA $L(T_0)$ using the following strategy.
The technique developed in [BB] provides a method to calculate any given homogeneous component
of the kernel of the epimorphism $\pi: \ V_\sg (c) \rightarrow L(T_0)$, though of course such a computation
is feasible only for the components of low degrees. Any $v \in \Ker \pi$ yields a relation between 
the fields in $L(T_0)$:
$$ Y_{L(T_0)} (v, z) = 0.$$
It turns out that it is sufficient to know the elements of $\Ker (\pi)$ of degrees 1 and 2 in order
to completely determine the structure of $L(T_0)$ as a VOA.

 Let us illustrate the technique of [BB] with the following example. Fix $m\in \Z^N$ and $g \in \dg$. 
Consider the subspace
$ P_1 = \Span \left< (t_0^{-1} t^r g) q^{m-r} \; | \; r\in\Z^N \right> \subset V_\sg(c) $. 
This subspace belongs to the
homogeneous component of weight $(-1,m)$ in $V_\sg(c)$. 
We are going to find the intersection of $P_1$ with $\Ker \pi$.
We note that a vector $v$ in component $(-1, m)$ of $V_\sg (c)$
belongs to $\Ker \pi$ if and only if $U_1(\g_+) v = 0$.

{\bf Lemma \LAA}. 
{\it
Let $g,g^\prime\in \dg$, $b = 1, \ldots, N$. Then
$$ (t_0 t^s k_0) (t_0^{-1} t^r g) q^{m-r}  = 0, \quad
 (t_0 t^s k_b) (t_0^{-1} t^r g) q^{m-r} = 0,  \eqno{(\kbg)} $$
$$ (t_0 t^s g^\prime) (t_0^{-1} t^r g) q^{m-r} = (g^\prime |g) c q^{m+s},  \eqno{(\gpg)} $$
$$ (t_0 t^s \td_0) (t_0^{-1} t^r g) q^{m-r} = 0,  \quad
 (t_0 t^s \td_b) (t_0^{-1} t^r g) q^{m-r} = 0.  \eqno{(\kdbg)} $$
}

{\it Proof.} Let us prove (\gpg):
$$ (t_0 t^s g^\prime) (t_0^{-1} t^r g) q^{m-r} = [t_0 t^s g^\prime, t_0^{-1} t^r g] q^{m-r}
+  (t_0^{-1} t^r g) (t_0 t^s g^\prime) q^{m-r} .$$
The second term vanishes because $\g_+$ acts trivially on the top of $V_\sg(c)$. 
For the first term we use the relations in $\g$,
(\koq) and (\kpq):
$$\eqalign{ [t_0 t^s g^\prime, t_0^{-1} t^r g] q^{m-r}
& = (t^{r+s} [g^\prime, g]) q^{m-r} + 
(g^\prime |g) (t^{r+s} k_0) q^{m-r} + (g^\prime |g)
\sum_{p=1}^N s_p (t^{r+s} k_p) q^{m-r} \cr
& = (g^\prime |g) c q^{m+s} \cr} .$$
All other equalities in the statement of this Lemma are obtained in a similar way.

 Let us analyze the results of this Lemma. We see that the right hand sides in (\kbg)-(\kdbg) are
independent of $r$. Thus 
$U_1(\g_+) \left( (t_0^{-1} t^r g) q^{m-r} - (t_0^{-1} g) q^m \right) = 0$,
which implies that $(t_0^{-1} t^r g) q^{m-r} - (t_0^{-1} g) q^m  \in \Ker \pi$.
Applying the state-field correspondence, we get that the following relation holds in $L(T_0)$:
$$ : Y(t_0^{-1} t^r g, z) Y(q^{m-r}, z) : \; = \; : Y(t_0^{-1} g, z) Y(q^m, z) : . \eqno{(\Ykg)}$$
Since these vertex operators commute, we may drop the normal ordering symbol.  
In particular, for $m=r$, we get that
$$   Y(t_0^{-1} t^m g, z)  =  Y(t_0^{-1} g, z) Y(q^m, z)  . \eqno{(\Ykkg)}$$
This factorization property means that the fields $Y(t_0^{-1} t^m g, z) = g(m,z)$ reduce to more elementary
fields $ Y(t_0^{-1} g, z) = g(z)$ and $Y(q^m, z)$. The fields $ g(z)$ correspond to the affine 
subalgebra $\C[t_0, t_0^{-1}] \otimes \dg \oplus \C k_0 \subset \g$. 
Below we establish analogous factorization formulas for other fields.
We shall also see that the affine fields $Y(t_0^{-1} g, z)$ commute with other elementary fields,
except for $Y(t_0^{-2} \td_0,z)$. This will
imply that the affine VOA generated by the affine fields splits off as a tensor factor in $L(T_0)$
when $c$ is not the critical level for this affine subalgebra. In this way
we will obtain a tensor product decomposition of $L(T_0)$.

Now let us get the formula for the vertex operators $Y(t_0^{-1} t^r k_a, z)$ and $Y(t^r k_0, z)$.
 For a fixed $m\in\Z^N$ and $1 \leq a \leq N$, consider the subspace
$ P_2 = \Span \left< (t_0^{-1} t^r k_a) q^{m-r}  | r\in\Z^N \right> \subset V_\sg(c) $.
Again we are going to find the intersection of $P_2$ with $\Ker \pi$. 

{\bf Lemma \LAB}. 
{\it
Let $1\leq a, b \leq N$. Then
$$ (t_0 t^s k_0) (t_0^{-1} t^r k_a) q^{m-r} = 0,  \quad
 (t_0 t^s k_b) (t_0^{-1} t^r k_a) q^{m-r}  = 0,  \eqno{(\kbao)} $$
$$ (t_0 t^s g) (t_0^{-1} t^r k_a) q^{m-r} = 0, \quad g \in \dg, \quad
 (t_0 t^s \td_0) (t_0^{-1} t^r k_a) q^{m-r}  = 0,  \eqno{(\kdoao)} $$
$$ (t_0 t^s \td_b) (t_0^{-1} t^r k_a) q^{m-r}  = c \delta_{ab} q^{m+s} .  \eqno{(\kdbao)} $$
}

 The proof of this Lemma is the same as for Lemma \LAA, and we omit these calculations.

\


 We see again that the right hand sides in (\kbao)-(\kdbao) are
independent of $r$. Thus 
$(t_0^{-1} t^r k_a) q^{m-r}  - (t_0^{-1} k_a) q^m  \in \Ker \pi$,
or equivalently,
$$ (t_0^{-1} t^r k_a) q^{m-r} = (t_0^{-1} k_a) q^m \quad {\rm in} \quad L(T_0). \eqno{(\kaq)}$$
Taking $r=m$ and 
applying the state-field correspondence, we get that the following relation holds in $L(T_0)$:
$$   Y(t_0^{-1} t^m k_a, z)  =  Y(t_0^{-1} k_a, z) Y(q^m , z)  . \eqno{(\Ykkao)}$$
Since these vertex operators commute, we dropped the normal ordering symbol in the right hand side.  
Also, taking into account that $t_0^{-1} t^m k_0 = \sum\limits_{p=1}^N m_p t_0^{-1} t^m k_p$, we obtain
$$  c \dz Y(q^m , z)  = Y( t_0^{-1} t^m k_0, z) =
 \sum_{p=1}^N m_p Y(t_0^{-1} k_p, z) Y(q^m, z)  . \eqno{(\Ykko)}$$

 Next we will use the results of Section 3 of [BB]. It is proved there that $L(T_0)$ is a module over a
commutative associative algebra $\C[q_1^\pm,\ldots,q_N^\pm]$ and the vertex operator 
$Y(q^m , z)$ is given by the expression (\Yko).
We can see that this formula is compatible with (\Ykko), and in fact it is not too difficult to derive (\Yko)
from (\Ykko). 

 We will later need another relation in $L(T_0)$ which can be derived either from (\kaq) or (\Yko):
$$ (t_0^{-1} t^r k_0) q^{s} = \sum_{p=1}^N r_p (t_0^{-1} k_p) q^{r+s} . \eqno{(\koqq)}$$

 Our next goal is to derive a formula for the vertex operator $Y(t_0^{-1} t^m \td_a, z)$. We will use the same
strategy as above.

{\bf Lemma \LAC.} 
{\it
Let $1\leq a, b \leq N$. Then
$$ (t_0 t^s k_0) (t_0^{-1} t^r \td_a) q^{m-r}  = -s_a c q^{m+s} ,  \eqno{(\koda)} $$
$$ (t_0 t^s k_b) (t_0^{-1} t^r \td_a) q^{m-r}  = \delta_{ab} c q^{m+s} ,  \eqno{(\kbda)} $$
$$ (t_0 t^s g) (t_0^{-1} t^r \td_a) q^{m-r} = 0, \quad g \in \dg, \eqno{(\gda)} $$
$$ (t_0 t^s \td_0) (t_0^{-1} t^r \td_a) q^{m-r}  = 
\left( (m_a - r_a) - 2 (\mu s_a - \nu r_a) c \right) q^{m+s} ,  \eqno{(\doda)} $$
$$ (t_0 t^s \td_b) (t_0^{-1} t^r \td_a) q^{m-r}  =
\left( r_b (m_a-r_a) - s_a (m_b - r_b) - (\mu r_b s_a + \nu r_a s_b) c \right) q^{m+s} .  \eqno{(\dbda)} $$
}

 {\it Proof.} Let us show the calculations for (\doda) and leave rest as
an exercise to the reader. We are going to use (\tdodb), (\koq) and (\tdpq):

$$ (t_0 t^s \td_0) (t_0^{-1} t^r \td_a) q^{m-r}  = [t_0 t^s \td_0, t_0^{-1} t^r \td_a] q^{m-r} $$ 
$$ = (t^{r+s} \td_a) q^{m-r}  - s_a (t^{r+s} \td_0) q^{m-r}  $$
$$+ ( -2 \mu s_a + 2 \nu r_a ) (t^{r+s} k_0) q^{m-r}   
- (- \mu s_a + 2 \nu r_a) \sum_{p=1}^N r_p (t^{r+s} k_p) q^{m-r}$$
$$= (m_a - r_a - 2\mu c s_a + 2\nu c r_a) q^{m+s}  .$$ 


Unlike the previous cases, the right hand sides in Lemma \LAC \ do depend on $r$. 
Note, however, that this dependence is polynomial, and we
may separate constant, linear and quadratic components. Clearly, the constant term is obtained by setting
$r=0$,
and is produced by the element $(t_0^{-1} \td_a) q^{m}$.
 
Next, comparing (\koda)-(\dbda) with (\kbao)-(\kdbao) we notice that the quadratic term is given by the vector
$-{r_a \over c} \sum\limits_{p=1}^N r_p (t_0^{-1} k_p) q^{m} $:
$$ (t_0 t^s \td_b) \left( -{r_a \over c} \sum_{p=1}^N r_p (t_0^{-1} k_p) q^{m}  \right) = 
-r_a r_b q^{m+s},  $$
while the other raising operators annihilate this vector.

Finally, the linear in $r$ component is given by the vector $\sum\limits_{p=1}^N r_p E_{pa}^m$, where
$$ E_{pa}^m = (t_0^{-1} t^{\epsilon_p} \td_a) q^{m-\epsilon_p} - (t_0^{-1} \td_a) q^{m} 
+ {1 \over c} \delta_{ap} (t_0^{-1} k_p) q^m . \eqno{(\Epam)}$$
Using Lemmas \LAC \ and \LAB, we can easily see that
$$ (t_0 t^s k_0) \sum_{p=1}^N r_p E_{pa}^m = 0, \quad
 (t_0 t^s k_b) \sum_{p=1}^N r_p E_{pa}^m = 0, \quad
 (t_0 t^s g) \sum_{p=1}^N r_p E_{pa}^m = 0, \eqno{(\gE)}$$
$$ (t_0 t^s \td_0) \sum_{p=1}^N r_p E_{pa}^m = r_a (-1 + 2\nu c) q^{m+s} , \eqno{(\doE)}$$
$$ (t_0 t^s \td_b) \sum_{p=1}^N r_p E_{pa}^m = 
( r_b m_a + (1 - \mu c) s_a r_b  - \nu c r_a s_b) q^{m+s} . \eqno{(\dbE)}$$

It follows from the above computations that the following vector is annihilated 
by $U_1(\g_+)$ and thus vanishes in $L(T_0)$:
$$ (t_0^{-1} t^r \td_a) q^{m-r} - (t_0^{-1} \td_a) q^{m} - \sum\limits_{p=1}^N r_p E_{pa}^m + 
{r_a \over c} \sum_{p=1}^N r_p (t_0^{-1} k_p) q^m  = 0 . \eqno{(\daq)}$$
 Setting $r=m$, and applying the map $Y$
we obtain a relation for the vertex operators in $L(T_0)$:
$$ \td_a (m,z) = :\td_a(z) Y(q^m,z): + \sum_{p=1}^N m_p Y( E_{pa}^m, z) 
- {m_a\over c} \sum_{p=1}^N m_p k_p(z) Y(q^m,z). \eqno{(\damz)}$$

 Let us now establish a relation between $E_{ab}^m$ and 
$$E_{ab} = E_{ab}^0 =  (t_0^{-1} t^{\epsilon_a} \td_b) q^{-\epsilon_a} - (t_0^{-1} \td_b) \o 
+ {1 \over c} \delta_{ab} (t_0^{-1} k_a) \o. \eqno{(\Eabo)}$$

{\bf Lemma \LAD.} 
{\it
(a) The fields $Y(q^m,z), g(z), k_p(z), \td_p (z)$ commute with
$Y(E_{ab},z)$, $\; p,a,b = 1, \ldots, N$.

(b) The following relation holds: 
$$E_{ab}^m = (E_{ab})_{(-1)} q^m + {m_b \over c} (t_0^{-1} k_a) q^m. \eqno{(\Emab)}$$ 
}

 {\it Proof.} Let us show that 
$$(q^m)_{(n)} E_{ab} = 0 \hbox{\rm \ for all \ } n\geq 0. \eqno{(\qE)}$$ 
Since $\deg ((q^m)_{(n)} E_{ab}) = -n$, we only need to consider
the case of $n=0$. But $(q^m)_{(0)} = {1\over c} (t_0 t^m k_0)$, and we get the desired claim from (\gE).
 Now applying the commutator formula (\comm) we get that the fields $Y(q^m,z)$ and $Y(E_{ab},z)$ commute.
Using a similar argument we can derive from (\gE)-(\dbE) that $g(z), k_p(z), \td_p (z)$ also commute
with $Y(E_{ab},z)$. 
    
 Taking into account the skew symmetry identity (\skeww) we obtain as a consequence of (\qE)
the equality
$$ (E_{ab})_{(-1)} q^m = (q^m)_{(-1)} E_{ab} .$$
If we substitute (\Eabo) in the right hand side of this equality, we will get
$$ (E_{ab})_{(-1)} q^m = 
{1\over c} (t^m k_0) \left( (t_0^{-1} t^{\epsilon_a} \td_b) q^{-\epsilon_a} - (t_0^{-1} \td_b) \o 
+ {1 \over c} \delta_{ab} (t_0^{-1} k_a) \o \right)$$
$$ = {1\over c} (t_0^{-1} t^{\epsilon_a} \td_b) (t^m k_0) q^{-\epsilon_a} - {1\over c} (t_0^{-1} \td_b) (t^m k_0) \o
+ {1 \over c^2} \delta_{ab} (t_0^{-1} k_a) (t^m k_0) \o$$
$$- {1\over c} [t_0^{-1} t^{\epsilon_a} \td_b, t^m k_0]  q^{-\epsilon_a} + {1\over c}  [t_0^{-1} \td_b, t^m k_0] \o$$
$$ = (t_0^{-1} t^{\epsilon_a} \td_b) q^{m-\epsilon_a} - (t_0^{-1} \td_b) q^m 
+ {1 \over c} \delta_{ab} (t_0^{-1} k_a) q^m - {m_b \over c} (t_0^{-1} t^{m+\epsilon_a} k_0) q^{-\epsilon_a}
+  {m_b \over c} (t_0^{-1} t^{m} k_0) \o$$
$$ = E_{ab}^m - {m_b \over c} (t_0^{-1} k_a) q^m .$$
To get the last equality we used (\Epam) and (\koqq). This completes the proof of the Lemma.

Combining (\damz) with (\Emab), we obtain (\Ydam). We also get from (\daq) and (\Emab) that
$$ (t_0^{-1} t^r \td_a) q^m = (t_0^{-1} \td_a) q^{m+r} + \sum_{p=1}^N r_p (E_{pa})_{(-1)} q^{m+r} + 
{m_a \over c}  \sum_{p=1}^N r_p (t_0^{-1} k_p) q^{m+r}. \eqno{(\daqq)} $$

Next we are going to determine the commutator relations between $Y(E_{ab},z)$, $1 \leq a,b \leq N$.

{\bf Lemma \LAE.} 
{\it
(a) $(E_{ab})_{(0)} E_{sp} = \delta_{bs} E_{ap} - \delta_{ap} E_{sb}$,

(b) $(E_{ab})_{(1)} E_{sp} = (1-\mu c) \delta_{bs} \delta_{ap} \o - \nu c  \delta_{ab} \delta_{sp} \o $, 

(c) $(E_{ab})_{(n)} E_{sp} = 0$ for $n \geq 2$.
}

{\it Proof.} Since $\deg(E_{ab}) = 1$, we have that $\deg\left( (E_{ab})_{(n)} E_{sp} \right) = 1 - n$,
from which part (c) immediately follows.
Let us carry out the computations for part (b) of the Lemma. We have
$$ (E_{ab})_{(1)} E_{sp} = \left( (t_0^{-1} t^{\epsilon_a} \td_b)_{(-1)} q^{-\epsilon_a} \right)_{(1)} E_{sp}
- (t_0 \td_b) E_{sp} + {1\over c} \delta_{ab} (t_0 k_a) E_{sp} . \eqno{(\EE)} $$
The last two terms vanish by (\gE) and (\dbE). To evaluate the first term in the right hand side
of (\EE), we use the Borcherds' identity (\Borb). Noting that by (\qE), $(q^{-\epsilon_a})_{(n)} E_{sp} = 0$
for all $n \geq 0$, we get
$$\left( (t_0^{-1} t^{\epsilon_a} \td_b)_{(-1)} q^{-\epsilon_a} \right)_{(1)} E_{sp} =
{1 \over c} (t_0 t^{-\epsilon_a} k_0) (t^{\epsilon_a} \td_b) E_{sp}
+ {1 \over c} (t^{-\epsilon_a} k_0) (t_0 t^{\epsilon_a} \td_b) E_{sp} .$$
The first term is equal to zero since
$$ {1 \over c} (t_0 t^{-\epsilon_a} k_0) (t^{\epsilon_a} \td_b) E_{sp} =
{1 \over c} (t^{\epsilon_a} \td_b) (t_0 t^{-\epsilon_a} k_0) E_{sp}
- {1 \over c} [ t^{\epsilon_a} \td_b, t_0 t^{-\epsilon_a} k_0 ] E_{sp} $$
and we may apply (\gE) to the right hand side. 
Finally, using (\dbE) and (\koq), we get
$${1 \over c} (t^{-\epsilon_a} k_0) (t_0 t^{\epsilon_a} \td_b) E_{sp} =
\left( (1-\mu c) \delta_{bs} \delta_{ap} - \nu c  \delta_{ab} \delta_{sp} \right) \o.$$
This completes the proof of part (b) of the Lemma, and we leave part (a) as an exercise for the reader.

Comparing Lemmas {\LAE} and {\reln} (a), we conclude that the operators $(E_{ab})_{(n)}$ produce 
a representation of affine $\wgl$.


{\bf Lemma \LAF.} 
{\it
The following relations hold in $L(T_0)$:
$$ (t_0^{-1} t^r \td_0) q^m = {1\over c} \sum_{p=1}^N m_p (t_0^{-1} k_p) q^{m+r} , \leqno{(a)}$$
$$ (t_0 t^r \td_0) (E_{ab})_{(-1)} q^m = \delta_{ab} (-1 + 2\nu c) . \leqno{(b)} $$
}

{\it Proof.}
$$ (t_0^{-1} t^r \td_0) q^m = {1\over c} (t_0^{-1} t^r \td_0) (t^m k_0) \o
= {1\over c} [ t_0^{-1} t^r \td_0, t^m k_0] \o $$
$$=  {1\over c} \sum_{p=1}^N m_p (t_0^{-1} t^{r+m} k_p) \o = {1\over c} \sum_{p=1}^N m_p (t_0^{-1} k_p) q^{r+m} .$$
This proves the claim (a). Part (b) follows from Lemma \LAD (b), (\doE) and (\kdoao).

Finally let us study the properties of the field $\td_0 (r,z)$.
This will require carrying out calculations with certain elements of degree $2$. Let $v \in L(T_0), \deg(v) =2$.
Since $L(T_0)$ is an irreducible $\g$-module, it is generated by any non-zero vector. If $v \neq 0$ then
$U_2 (\g_+) v = T_0 $.
However it is easy to see that $\g_+$ is generated by
$\g_1$. Thus $\g_1 v = 0$ implies that $v = 0$ . This observation will help us find a relation in
$L(T_0)$ involving $(t_0^{-2} t^m \td_0) \o$, which in turn will yield an expression for the field $\td_0 (m,z)$.

{\bf Lemma \LAG.} 
{\it
The following relation holds in $L(T_0)$:
$$ (t_0^{-2} t^m \td_0) \o = (t_0^{-2} \td_0) q^m + {1\over c} \sum_{p,j=1}^N m_p (t_0^{-1} k_j)
(E_{pj})_{(-1)} q^m
- {1\over c} (1-\mu c) \sum_{p=1}^N m_p (t_0^{-2} k_p) q^m .$$ 
}

{\it Proof.} We shall prove this Lemma by showing that the vector 
$$ v = (t_0^{-2} t^m \td_0) \o - (t_0^{-2} \td_0) q^m - {1\over c} \sum_{p,j=1}^N m_p (t_0^{-1} k_j)
(E_{pj})_{(-1)} q^m
+ {1\over c} (1-\mu c) \sum_{p=1}^N m_p (t_0^{-2} k_p) q^m $$ 
is annihilated in $L(T_0)$ by $\g_1$. Let us show that $(t_0 t^s \td_0) v = 0$ in $L(T_0)$: 
$$ (t_0 t^s \td_0) v = [ t_0 t^s \td_0, t_0^{-2} t^m \td_0] \o -  [ t_0 t^s \td_0, t_0^{-2} \td_0] q^m$$
$$ - {1\over c} \sum_{p,j=1}^N m_p [ t_0 t^s \td_0, t_0^{-1} k_j ] (E_{pj})_{(-1)} q^m
- {1\over c} \sum_{p,j=1}^N m_p (t_0^{-1} k_j) (t_0 t^s \td_0) (E_{pj})_{(-1)} q^m$$
$$ + {1\over c} (1-\mu c) \sum_{p=1}^N m_p [t_0 t^s \td_0, t_0^{-2} k_p] q^m $$
$$ = 3 (t_0^{-1} t^{m+s} \td_0) \o + 4(\mu+\nu) (t_0^{-1} t^{m+s} k_0)\o 
- 2(\mu+\nu) \sum_{p=1}^N m_p (t_0^{-1} t^{m+s} k_p) \o$$
$$- 3 (t_0^{-1} t^{s} \td_0) q^m - 4(\mu+\nu) (t_0^{-1} t^{s} k_0) q^m $$
$$- {1\over c} \sum_{p,j=1}^N m_p (t^s  k_j) (E_{pj})_{(-1)} q^m
- {1\over c} (-1 + 2\nu c) \sum_{p=1}^N m_p (t_0^{-1} k_p) q^{m+s}$$
$$+ {2 \over c} (1-\mu c)  \sum_{p=1}^N m_p  (t_0^{-1} t^s k_p) q^{m}$$
$$ = 4(\mu+\nu) \sum_{p=1}^N (m_p + s_p) (t_0^{-1} k_p) q^{m+s}
- 2(\mu+\nu) \sum_{p=1}^N m_p (t_0^{-1} k_p) q^{m+s} - {3\over c} \sum_{p=1}^N m_p (t_0^{-1} k_p) q^{m+s} $$
$$ - 4(\mu+\nu) \sum_{p=1}^N s_p (t_0^{-1} k_p) q^{m+s} 
- {1\over c} \sum_{p,j=1}^N m_p  (E_{pj})_{(-1)} (t^s  k_j) q^m$$
$$ - {1\over c} (-1 + 2\nu c) \sum_{p,j=1}^N m_p (t_0^{-1} k_p) q^{m+s}
+ {2 \over c} (1-\mu c)  \sum_{p=1}^N m_p  (t_0^{-1} k_p) q^{m+s} = 0. $$
The cases of other elements
spanning $\g_1$ are treated similarly, and are left as an exercise.

Applying the state-field correspondence $Y$ to both sides of the equality in Lemma \LAG, we will get the formula (\Ydom).  

Let us complete the proof of Theorem \TA. 
We have now established all the relations between the fields stated in part (c) of the Theorem. The universal enveloping
vertex algebra $V_\sg$ is generated by the fields (\korz)-(\dorz). The same is true for  $L(T_0)$,
since $L(T_0)$ is a factor vertex algebra of $V_\sg$. Taking into account the relations of part (c), 
the claim of part (a)
of the Theorem follows. 
As we mentioned above, the claim of part (b) follows from the results of [BB].

Let us establish the claim of part (d) of the Theorem. 
First we construct a homomorphism of vertex algebras
$$ \varphi: \quad \VH \rightarrow L(T_0), $$
defined by $\varphi (e^{mu}) = q^m, \; \varphi(u_p(-1) \o) = {1 \over c} (t_0^{-1} k_p) \o, \; 
\varphi(v_p(-1)\o) = (t_0^{-1} \td_p) \o $. 
Using (\Ldk), (\kpq), (\dpq) and (\Yko), we can see that the images of the generators of
$\VH$ satisfy the required relations (\xyK)-(\Ye), thus
$\varphi$ is indeed a homomorphism of vertex algebras. Since $\VH$ is a 
simple vertex algebra, the map $\varphi$ is injective. The image of the Virasoro field of $\VH$
is
$$ \varphi(\omega_\Hyp (z)) = {1\over c} \sum_{p=1}^N :d_p(z) k_p(z): ,$$
and the central charge of this Virasoro field is equal to $\rank(\VH) = 2N$.

We know that the fields $g(z)$ generate an affine subalgebra $\wdg$ in the toroidal algebra $\g$,
which commutes with the fields generated by the image of $\varphi$. This affine subalgebra
has the central charge $c_\sdg = c$ on $L(T_0)$.
Next, comparing Lemmas {\LAE}  and \reln (a), we conclude that the fields 
$Y(E_{ab},z)$ yield a representation of the affine $\wgl$ on $L(T_0)$ with central charges
$c_\sln = 1 - \mu c$ and $c_\Hei = N(1-\mu c) - N^2 \nu c$.
It follows from Lemma {\LAD} (a) that the fields $Y(E_{ab},z)$ also commute with the image of $\varphi$.

The relation (\tdodo) implies that the field $\td_0(z)$ generates  
a Virasoro algebra with the central charge $12 (\mu + \nu) c$. 
The fact that the element $t_0^{-1} \td_0$
is an infinitesimal translation operator $D$ follows from the relations (\Ldg), (\Ldk), (\tdodb) and (\tdodo).
Thus $L(T_0)$ is a VOA of rank $12 (\mu + \nu) c$ with the Virasoro field $\td_0(z)$.
The field $\td_0(z) - \varphi(\omega_\Hyp (z))$ yields another Virasoro algebra with central
charge $12 (\mu + \nu) c - 2N$. Comparing Lemmas {\reln} (b) and {\LAF} (b), and noting that
$\varphi(\omega_\Hyp (z))$ commutes with $g(z)$ and $Y(E_{ab},z)$, we get that the fields
$\td_0(z) - \varphi(\omega_\Hyp (z))$, $g(z)$ and $Y(E_{ab},z)$ yield a representation of the
twisted Virasoro-affine algebra $\f$ with $\df = \dg \oplus \glN$ and central character
given by $(\gamo)$. 


 This allows us to define a homomorphism of vertex algebras
$$ \psi: V_\sf (\gamma_0) \rightarrow L(T_0)$$
by $\psi(g(-1) \o) = (t_0^{-1} g) \o$, $\psi(E_{ab}(-1) \o) = E_{ab}$,
$\psi(\omega_\sf) = (t_0^{-2} \td_0) \o - \varphi(\omega_\Hyp)$. 
The sub-VOAs
$\varphi(\VH)$ and $\psi(V_\sf)$ commute in $L(T_0)$. Thus we have a homomorphism
$$\theta: \quad \VH \otimes V_\sf (\gamma_0) \rightarrow L(T_0) .$$
Since $\theta( \omega_\Hyp + \omega_\sf) = (t_0^{-2} \td_0) \o$, this is in fact a homomorphism 
of VOAs. Moreover, by part (a), $\theta$ is an epimorphism.

The $\f$-module $V_\sf (\gamma_0)$ has a unique maximal submodule and a unique 
irreducible quotient $L_\sf(\gamma_0)$. By Theorem {\inv}, $L_\sf(\gamma_0)$ is a unique quotient 
vertex algebra of 
$V_\sf (\gamma_0)$. Since $\VH$ is a simple vertex algebra, we conclude that
$\VH \otimes L_\sf(\gamma_0)$ is a unique simple quotient vertex algebra of 
$ \VH \otimes V_\sf (\gamma_0)$.
However $L(T_0)$ is also a simple quotient of $ \VH \otimes V_\sf (\gamma_0)$.
Thus
$$ L(T_0) \cong \VH \otimes L_\sf(\gamma_0). $$
This completes the proof of Theorem \TA.

\

{\bf 5. Realizations for the irreducible modules in category $\B_\chi$.}

\
 
 In this section we are going to use the theory of VOA modules to give realizations for all 
irreducible $\g(\mu,\nu)$-modules in category $\B_\chi$. By the principle of preservation of 
identities [Li], every VOA module for $ \VH \otimes L_\sf(\gamma_0)$ is also a module for 
the Lie algebra $\g(\mu,\nu)$. However in order to get all irreducible modules in $\B_\chi$,
we need to use a larger VOA 
$$V(T_0) = \VH \otimes V_\sf (\gamma_0) ,$$ 
which  has more irreducible modules than  $ \VH \otimes L_\sf(\gamma_0)$. 
In order to carry out this plan, we should first prove that $\VH \otimes V_\sf (\gamma_0)$ 
also admits a structure of a $\g(\mu,\nu)$-module.

First, we need to establish the following technical lemma:

{\bf Lemma \PAB.} 
{\it
For a Zariski dense set of triples $(c,\mu,\nu)$,
the modules $V_\sf (\gamma_0)$ and $L_\sf (\gamma_0)$ coincide.
}

{\it Proof.} It follows from Proposition {\FA} and (\gamo) that whenever 
$c \neq 0$, $c \neq -h^\vee$, $c_\sln = 1-\mu c \neq -N$, 
$c_\Hei = N(1-\mu c) - N^2 \nu c \neq 0$,
the VOA $V_\sf (\gamma_0)$ factors into a tensor product of four VOAs:
$V_\swdg (c_\sdg)$, $V_\wsl (c_\sln)$, $V_\Hei (c_\Hei)$ and $V_\Vir (c^\prime_\Vir)$.
It is clear that $V_\sf (\gamma_0)$ is a simple VOA if and only if each of these four
VOAs is simple. First of all, under the assumption that $c_\Hei \neq 0$,
the Heisenberg VOA is simple. We are going to show that the remaining affine and 
Virasoro VOAs are simple for $(c, \mu, \nu)$ in a dense subset of $\C^3$.

We note that affine and Virasoro VOAs are the generalized Verma modules for the 
respective Lie algebras, and a generalized Verma module admits a Shapovalov form [J].
Using the Shapovalov determinant argument, it is not difficult to see that
the VOA  $V_\swdg (c_\sdg)$ (resp. $V_\wsl (c_\sln)$, $V_\Vir (c^\prime_\Vir)$)
is simple outside a countable set of values of the central charge $c_\sdg$
(resp. $c_\sln$, $c_\Vir^\prime)$. 
In fact, an explicit formula for the Shapovalov determinant for the generalized Verma modules
for the affine algebras may be found in [KK], while for the Virasoro VOA, it follows from the 
description of the irreducible modules for the Virasoro algebra [FF] 
that $V_\Vir (c_\Vir^\prime)$ is simple if and 
only if $c_\Vir^\prime \neq 1 - 6{ (r-s)^2 \over rs}$, where $r$ and $s$ are relatively prime
integers, $r,s > 1$. 

Each inequality on the values of the central charges defines a dense Zariski open subset
of values of $(c, \mu, \nu)$. Since a countable intersection of dense Zariski open subsets
in $\C^3$ is Zariski dense, we establish the claim of the Lemma.

\

{\bf Proposition \PA.} Let $c\neq 0$, and let $\gamma_0$ be given by (\gamo).
Then $V(T_0) = \VH \otimes V_\sf(\gamma_0)$ has a structure of a $\g(\mu,\nu)$-module given by 
the formulas of Theorem {\TA} (b)-(d).

{\it Proof.} We proved in Theorem {\TA} (d) that $\VH \otimes L_\sf (\gamma_0)$ is 
a $\g(\mu,\nu)$-module. Now we
want to prove the same for  $\VH \otimes V_\sf(\gamma_0)$. This amounts to verifying relations
(\Rkakb)-(\Rdodo) in this vertex algebra. 
It is possible to do this directly, and this was the approach taken, for example, in [B4], 
but we are going to present here an alternative argument that allows us 
to circumvent these rather tedious computations.

 By Lemma \PAB, for a Zariski dense set of triples $(c,\mu,\nu)$,
the modules $V_\sf (\gamma_0)$ and $L_\sf (\gamma_0)$ coincide. Thus for the generic values of 
$(c,\mu,\nu)$
the VOAs $\VH \otimes V_\sf(\gamma_0)$ are indeed $\g(\mu,\nu)$-module, and the relations 
(\Rkakb)-(\Rdodo) hold.
However the commutator formula (\comm) applied to the left hand sides of the relations 
(\Rkakb)-(\Rdodo)
in $\VH \otimes V_\sf(\gamma_0)$ 
will yield expressions with coefficients that
are polynomials in $c^{\pm 1}, \mu, \nu$. Since these agree with the right hand sides of 
(\Rkakb)-(\Rdodo)
on a Zariski dense set of parameters, the equalities must hold for all values of $(c,\mu,\nu)$ 
with $c\neq 0$.
This concludes the proof of the Proposition. 

 Now we are ready to give realizations for all irreducible $\g$-modules in category $\B_\chi$ 
using highest weight $\f$-modules.


{\bf Theorem \TC.} 
{\it Let $c\neq 0$, and let $L(T)$ be an irreducible module in category $\B_\chi$
determined by the data $(V,W,h,d,\alpha)$ as in Theorem \BC, where
$V$ is a  finite-dimensional irreducible $\dg$-module, 
$W$ is a finite-dimensional irreducible $\slN$-module,
$\alpha \in \C^N $ and $h, d\in \C$. 
Then
$$ L(T) \cong M_\Hyp^+ (\alpha) \otimes L_\sf (V, W, h_\Hei, h_\Vir, \gamma_0), $$
where $\gamma_0$ is the same as in Theorem \TA,
$$ h_\Hei = h - N\nu c, \quad h_\Vir = -d + {1 \over 2} (\mu + \nu) c . \eqno{(\hhh)}$$
}

{\it Proof.} First of all, we are going to show that 
$M_\Hyp^+ (\alpha) \otimes L_\sf (V, W, h_\Hei, h_\Vir, \gamma_0)$
is a  $\g(\mu,\nu)$-module. Indeed, it is a module for the vertex algebra 
$\VH \otimes V_\sf(\gamma_0)$,
and by principle of preservation of identities [Li], the $\g(\mu,\nu)$-module structure on 
$\VH \otimes V_\sf(\gamma_0)$ gets transferred to its VOA module 
$M_\Hyp^+ (\alpha) \otimes L_\sf (V, W, h_\Hei, h_\Vir, \gamma_0)$.

Next, let us show that 
$ M_\Hyp^+ (\alpha) \otimes L_\sf (V, W, h_\Hei, h_\Vir, \gamma_0)$ is irreducible as a $\g$-module. 
This is not difficult to see. The fields that define the $\g$-module structure on
$\VH \otimes V_\sf(\gamma_0)$ generate this VOA. Thus any $\g$-submodule in 
$M_\Hyp^+ (\alpha) \otimes L_\sf (V, W, h_\Hei, h_\Vir, \gamma_0)$ is also a VOA submodule. 
However as a VOA module, $ M_\Hyp^+ (\alpha) \otimes L_\sf (V, W, h_\Hei, h_\Vir, \gamma_0)$ 
is irreducible. Thus it is also irreducible as a $\g$-module.

It is also easy to see that the $\g(\mu,\nu)$-module 
$M_\Hyp^+ (\alpha) \otimes L_\sf (V, W, h_\Hei, h_\Vir, \gamma_0)$
belongs to the category $\B_\chi$, and its top is 
$$\C [q_1^\pm, \ldots, q_N^\pm] \otimes V \otimes W. $$ 
We can derive from (\Yko)-(\Ydom) that $\g_0$ acts on this top according to (\kqu), (\dvw)-(\gvw).
\break
Taking into account the relations (\tdp) and (\tdo), we conclude that the top of
\break
$M_\Hyp^+ (\alpha) \otimes L_\sf (V, W, h_\Hei, h_\Vir, \gamma_0)$
is isomorphic to $T$ as a $\g_0$-module. Since two simple modules with the same top are
isomorphic, we obtain the claim of the Theorem.


Finally, applying Corollary \FAA, we obtain the following result:

{\bf Theorem \FB.} 
{\it
 Let $L(T)$ be the irreducible $\g(\mu,\nu)$-module in category $\B_\chi$ 
determined by the data $(V, W, h, d, \alpha)$ as in the statement of Theorem \BC.
Suppose
$$c \neq 0, \quad c \neq -h^{\vee}, \quad c_\sln = 1 -\mu c \neq -N, $$
$$c_\Hei = N(1-\mu c) - N^2 \nu c \neq 0,$$
$$c_\Vir^\prime = 12 c (\mu + \nu) - 2N - {c \dim(\dg) \over c + h^\vee}
- { c_\sln (N^2 -1) \over c_\sln + N} -1 + 12 {N^2 ({1\over 2} - \nu c)^2 \over c_\Hei} .$$
Let 
$$ h_\Hei = h - N\nu c, $$
$$ h_\Vir^\prime = -d + {1\over 2} (\mu + \nu) c - {\Omega_V \over 2(c + h^\vee)}
- {\Omega_W \over 2(c_\sln + N)}
- {h_\Hei (h_\Hei - N(1 - 2\nu c)) \over 2c_\Hei} .$$
Then
$$ L(T) \cong M_\Hyp (\alpha) \otimes L_{\swdg} (V, c) \otimes L_{\wsl} (W, c_\sln) 
\otimes L_{\Hei} (h_\Hei, c_\Hei) \otimes L_{\Vir} (h_\Vir^\prime, c_\Vir^\prime), \leqno{(a)}$$
$$ \eqalign {(b)  \quad
\char L(T) = &  
 \char q^\alpha \C[q_1^\pm,\ldots, q_N^\pm] \times \prod_{j \geq 1} (1-t^j)^{-(2N+1)} \cr
& \times
\char L_{\swdg} (V, c) \times \char L_{\wsl} (W, c_\sln) 
\times \char L_{\Vir} (h_\Vir^\prime, c_\Vir^\prime) .}$$
}

\

{\bf Remark \sll.} In case of 2-toroidal Lie algebras ($N=1$), the $\slN$ piece will
not be present, and should be omitted from all the statements in this paper.

\

\

{\bf References:}

\item{[ABFP]} B.~Allison, S.~Berman, J.~Faulkner and A.~Pianzola, 
{\it Realizations of graded-simple algebras as loop algebras,}
to appear.


\item{[ACKP]} E.~Arbarello, C.~De Concini, V.G.~Kac and C.~Procesi, 
{\it Moduli spaces of curves and representation theory,}
Comm.Math.Phys. {\bf 117}, 1-36 (1988).

\item{[BB]} S.~Berman and Y.~Billig, 
{\it Irreducible representations for toroidal Lie algebras,} 
J.Algebra {\bf 221}, 188-231 (1999).

\item{[BBS]} S.~Berman,Y.~Billig and J.~Szmigielski, 
{\it Vertex operator algebras and the representation theory of
toroidal algebras,}
in ``Recent developments in infinite-dimensional Lie algebras and conformal
field theory'' (Charlottesville, VA, 2000), 
Contemp. Math. {\bf 297}, 1-26, Amer. Math. Soc., 2002.

\item{[BC]} S.~Berman and B.~Cox,  
{\it  Enveloping algebras and representations of toroidal Lie algebras,} 
Pacific J.Math. {\bf 165}, 239-267 (1994).


\item{[B1]} Y.~Billig,  
{\it Principal vertex operator representations for toroidal Lie algebras,}  
J. Math. Phys. {\bf 39}, 3844-3864 (1998).

\item{[B2]} Y.~Billig,  
{\it An extension of the KdV hierarchy arising from a representation 
of a toroidal Lie algebra,} 
J.Algebra {\bf 217}, 40-64 (1999).

\item{[B3]} Y.~Billig, 
{\it Representations of the twisted Heisenberg-Virasoro algebra at level 
zero,}
Canadian Math. Bull. {\bf 46}, 529-537 (2003).

\item{[B4]} Y.~Billig, 
{\it Energy-momentum tensor for the toroidal Lie algebras,}
math.RT/0201313.

\item{[B5]} Y.~Billig, 
{\it Jet modules,} 
math.RT/0412119, to appear in Canadian J. of Math.

\item{[B6]} Y.~Billig, 
{\it Representations of toroidal extended affine Lie algebras,}
math.RT/0602112, to appear in J.Algebra.

\item{[BL]} Y.~Billig and M.~Lau, 
{\it Irreducible modules for extended affine Lie algebras,}
in preparation.

\item{[BN]} Y.~Billig and K.-H.~Neeb, 
{\it On the cohomology of vector fields on tori,}
in preparation.

\item{[DFP]} I.~Dimitrov, V.~Futorny and I.~Penkov,
{\it A reduction theorem for highest weight modules over toroidal Lie 
algebras,}
Comm. Math. Phys. {\bf 250}, 47-63 (2004).

\item{[DLM]} C.~Dong, H.~Li and G.~Mason, 
{\it Vertex Lie algebras, vertex Poisson algebras and vertex algebras,}
in ``Recent developments in infinite-dimensional Lie algebras and conformal
field theory'' (Charlottesville, VA, 2000), 
Contemp. Math. {\bf 297}, 69-96, Amer. Math. Soc., 2002.

\item{[E]} S.~Eswara Rao, 
{\it Partial classification of modules for Lie algebra of diffeomorphisms of
$d$-dimensional torus,}
J. Math. Phys. {\bf 45}, 3322-3333 (2004).

\item{[EM]} S.~Eswara Rao and R.V.~Moody, 
{\it Vertex representations for $n$-toroidal
Lie algebras and a generalization of the Virasoro algebra,}
Comm.Math.Phys. {\bf 159}, 239-264 (1994).

\item{[FF]} B.L.~Feigin and D.B.~Fuks, 
{\it Verma modules over the Virasoro algebra.}
Funktsional. Anal. i Prilozhen. {\bf 17}, 91-92 (1983).

\item{[FKRW]} E.~Frenkel, V.~Kac, A.~Radul and W.~Wang, 
{\it $W_{1+\infty}$ and $W(gl_\infty)$ with central charge $N$,}
Comm.Math.Phys. {\bf 170}, 337-357 (1995). 

\item{[FJW]} I.B.~Frenkel, N.~Jing and W.~Wang, 
{\it Vertex representations via finite groups and
the McKay correspondence,}
Internat. Math. Res. Notices {\bf 4}, 195-222 (2000).

\item{[FLM]} I.B.~Frenkel, J.~Lepowsky and A.~Meurman, 
{\it Vertex operator algebras and the Monster,} 
New York, Academic Press, 1988.

\item{[F]} D.B.~Fuks, 
{\it Cohomology of infinite-dimensional Lie algebras,}
New York, Consultants Bureau, 1986. 

\item{[IT]} T.~Ikeda and K.~Takasaki, 
{\it Toroidal Lie algebras and Bogoyavlensky's 2+1-dimensional equation,}
Internat.Math.Res.Notices {\bf 7}, 329-369 (2001).

\item{[IKU]} T.~Inami, H.~Kanno and T.~Ueno, 
{\it Higher-dimensional WZW model on K\"ahler manifold and toroidal Lie algebra,}
Mod.Phys.Lett. A {\bf 12}, 2757-2764 (1997).

\item{[IKUX]} T.~Inami, H.~Kanno, T.~Ueno and C.-S.~Xiong, 
{\it Two-toroidal Lie algebra as current algebra of four-dimensional K\"ahler WZW model,}
Phys.Lett. B {\bf 399}, 97-104 (1997).

\item{[ISW]} K.~Iohara, Y.~Saito and M.~Wakimoto,  
{\it Hirota bilinear forms with 2-toroidal symmetry,} 
Phys.Lett. A {\bf 254}, 37-46 (1999).

\item{[J]} J.C.~Jantzen, 
{\it Kontravariante Formen auf induzierten Darstellungen halbeinfacher
Lie-Algebren,}
Math. Ann. {\bf 226}, 53-65 (1977).

\item{[JM]} C.~Jiang and D.~Meng,
{\it Integrable representations for generalized Virasoro-toroidal Lie algebras,}
J. Algebra {\bf 270}, 307-334 (2003).

\item{[K1]} V.G.~Kac, 
{\it Infinite dimensional Lie algebras,} 
Cambridge, Cambridge University Press, 3rd edition, 1990.

\item{[K2]}  V.G.~Kac,  
{\it Vertex algebras for beginners,} 
2nd edition, University Lecture Series, {\bf 10}, Amer. Math. Soc., 1998.

\item{[KK]} V.G.~Kac and D.A.~Kazhdan, 
{\it Structure of representations with highest weight of
infinite-dimensional Lie algebras,}
Adv. Math. {\bf 34}, 97-108 (1979).

\item{[Kas]} C.~Kassel, 
{\it K\"ahler differentials and coverings of complex simple 
Lie algebras extended over a commutative ring,}
J.Pure Applied Algebra {\bf 34}, 265-275 (1984).

\item{[L]} T.A.~Larsson, 
{\it Lowest-energy representations of non-centrally extended diffeomorphism 
algebras,}
Comm.Math.Phys. {\bf 201}, 461-470 (1999).

\item{[Li]} H.~Li, 
{\it Local systems of vertex operators, vertex superalgebras and modules,}
{J.Pure Appl. Algebra} {\bf 109}, 143-195 (1996).

\item{[MRY]} R.V.~Moody, S.E.~Rao and T.~Yokonuma,   
{\it Toroidal Lie algebras and vertex representations,} 
Geom.Ded. {\bf 35}, 283-307 (1990).

\item{[P]} M.~Primc,
{\it Vertex algebras generated by Lie algebras,}
J.Pure Appl.Algebra {\bf 135}, 253-293 (1999).

\item{[R]} M.~Roitman, 
{\it On free conformal and vertex algebras,}
J.Algebra {\bf 217}, 496-527 (1999). 

\item{[T]} T.~Tsujishita, 
{\it Continuous cohomology of the Lie algebra of
vector fields,}
Mem.Amer. Math.Soc. {\bf 34} no.253 (1981).

\

\

{School of Mathematics \& Statistics,}
{Carleton University,}
{1125 Colonel By Drive,}
{Ottawa, Ontario, K1S 5B6, Canada}

{E-mail address: billig@math.carleton.ca}

\end